\documentclass[a4paper,10pt]{article}
\usepackage{mathrsfs}
\usepackage{amsmath}
\usepackage{amsfonts}
\usepackage{graphicx}
\usepackage{float}
\usepackage{amssymb}
\usepackage{geometry}
\usepackage{cases}
\usepackage{enumerate}
\usepackage{setspace}

\usepackage{color}
\usepackage{arydshln}

\newtheorem{example}{Example}[section]
\newtheorem{proposition}{Proposition}[section]

\newtheorem{definition}{Definition}[section]
\newtheorem{theorem}{Theorem}[section]
\newtheorem{lemma}{Lemma}[section]
\newtheorem{remark}{Remark}[section]

\usepackage{color,xcolor}

\date{}

\begin{document}
\title{Stochastic Singular Linear Systems and Related Linear-Quadratic Optimal Control Problems under Finite and Infinite Horizons
\author{Mengzhen Li\footnote{{\it E-mail address}: 13335291065@163.com (M. Li).} \hspace{1cm} Tianyang Nie\footnote{{\it E-mail address}: nietianyang@sdu.edu.cn (T. Nie).}
\hspace{1cm} Zhen Wu\footnote{Corresponding author.
{\it E-mail address}: wuzhen@sdu.edu.cn (Z. Wu).}\\
\vspace{1mm}
\small\em School of Mathematics, Shandong University, Jinan {\rm 250100}, China}}

\maketitle

\vspace{2mm}\noindent\textbf{Abstract.}\quad
In this paper, we study the necessary and sufficient conditions for ensuring the well-posedness of the stochastic singular systems. Moreover, we investigate the stochastic singular linear-quadratic control problems, considering both finite and infinite horizons, and transform each of these problems into their corresponding normal linear-quadratic control problem. To guarantee the finiteness of the infinite-horizon linear-quadratic control problem, we establish the Popov-Belevitch-Hautus rank criterion for accessing the controllability of the stochastic system. Furthermore, we derive the feedback form of the optimal control. Finally, we provide solutions for illustrated examples of the stochastic singular linear-quadratic control problem in both finite and infinite horizons.

\vspace{2mm} \noindent\textbf{Keywords.}\quad
Singular system, well-posedness, linear-quadratic optimal control problem, exactly controllable, infinite horizon, maximum principle.

\vspace{2mm}\noindent\textbf{AMS subject classification.} 93E20, 60H10


\section{Introduction}
Since Rosenbrock \cite{rosenbrock1974structural} introduced the concept of singular systems while studying complex electrical circuit network, research in this field has flourished across various disciplines, including social systems \cite{kumar1995feedback}, power systems \cite{lewis1986survey}, biological systems \cite{liu2008passivity}, economic systems \cite{pantelous2010linear} and other areas \cite{karcanias2008structure}. Singular systems are also known as descriptor systems, differential-algebraic systems, and generalized state-space systems, as elaborated in references like \cite{campbell1982singular,
zhaolin1988optimal,
luenberger1977dynamic}.
Now we consider the following deterministic linear system with time-invariant coefficients
\begin{equation}\label{ode}
\begin{aligned}
E\dot{x}(t)=Ax(t)+Bu(t),\\
\end{aligned}
\end{equation}
where $x(\cdot)$ is the state vector, $u(\cdot)$ is the control input, and $E,\ A,\ B$ are constant matrices of suitable dimensions.
When matrix $E$ is invertible, system \eqref{ode} turns into
\begin{equation*}
\begin{aligned}
&\dot{x}(t) = E^{-1}[Ax(t)+Bu(t)].
\end{aligned}
\end{equation*}
This is a common linear system, extensively studied and applied in various contexts.
In this paper, we refer to such systems as normal systems. When matrix $E$ is singular, the corresponding system \eqref{ode} is denoted as a singular system.  
For the normal linear system, there always exists a unique solution for any control. Conversely, for the singular system, the existence of a unique solution cannot be guaranteed.
\begin{example}
Consider the controlled singular  differential-algebraic equation (DAE)
\begin{equation*}
\left\{
\begin{aligned}
&\begin{pmatrix}
1  &0  &0\\
0  &0  &0\\
0  &0  &0
\end{pmatrix}
\begin{pmatrix}
\dot{x}_1\\
\dot{x}_2\\
\dot{x}_3
\end{pmatrix}
=\begin{pmatrix}
1  &0  &0\\
0  &1  &0\\
0  &0  &0
\end{pmatrix}
\begin{pmatrix}
x_1\\
x_2\\
x_3
\end{pmatrix}
+\begin{pmatrix}
u_1\\
u_2\\
u_3
\end{pmatrix},\\
&\quad x_1(0)=1,\ x_2(0)=-u_2(0),\ x_3(0)=1.
\end{aligned}
\right.
\end{equation*}
Then, evidently, the DAE has a solution only when $u_3=0$.
\end{example}

Roughly speaking, a singular system admits a unique solution for any control only when specific conditions are met by the system coefficients. Additionally, due to the presence of a singular matrix on the left-hand side of the equation, Gronwall's inequality can no longer be used to estimate the trajectories. Consequently, compared to a normal system, a singular system faces at least two crucial obstacles. First, it is impossible to establish general existence and uniqueness results due to their complex structure. Second, the maximum principle fails for singular systems. Many efforts have been made to address the singular linear-quadratic (LQ) control problem in the deterministic case. Interested readers can refer to
\cite{bender1987linear,
zhaolin1988optimal,
cobb1983descriptor,
dai1989singular,
duan2010analysis,
zhu2002singular}. In these works,
the main method involves transforming the singular system into a normal system through equivalent transformation.

In actual models, random noise cannot be ignored. It can provide a more accurate description of real-world models by incorporating random noise into singular systems. For example, a portfolio of various insurance products and the investigation of the pricing process can be modeled using linear time-invariant discrete-time stochastic singular systems, see \cite{gashi2013linear,
pantelous2010linear}.
In the past decades, a considerable amount of
papers have been published concerned with the stochastic singular systems \cite{boukas2006stabilization,
gao2013observer,gashi2013linear,
ge2021impulse,
huang2010stability,
zhang2014stability,
zhang2017linear,zhang2015some}, which mainly focused on the stability (stabilization) and controllability.
For example, Boukas \cite{boukas2006stabilization} assumed the existence of the solution to the stochastic singular hybrid system and considered the stochastic stability. Huang and Mao \cite{huang2010stability} imposed two sufficient conditions to ensure the well-posedness of the stochastic singular system with Markovian switching.
Gao and Shi \cite{gao2013observer} pointed out that the regularity can ensure the existence of the solution to the singular stochastic differential equations (SDEs) but without any proof and addressed mean-square exponential stability for stochastic singular systems. Gashi and Pantelous \cite{gashi2013linear} established the necessary and sufficient conditions for ensuring the well-posedness of linear regular
backward stochastic differential equations (BSDEs). Zhang et al. \cite{zhang2015some} introduced two different sufficient conditions for guaranteeing the well-posedness of uncontrolled singular SDEs and dealt with the stability for stochastic singular systems in both continuous-time and discrete-time case. Ge \cite{ge2021impulse} gave the sufficient condition for ensuring the well-posedness of stochastic singular systems and discussed the topics related to impulse controllability and impulse observability.

However, as far as we know, there have been only two papers devoted to studying the stochastic singular optimal control problem in the past. Zhang and Xing \cite{zhang2014stability} provided sufficient conditions for the well-posedness of uncontrolled singular SDE and considered stochastic singular LQ control problem under both finite and infinite horizons, respectively. Zhang et al. \cite{zhang2017linear} offered sufficient condition for the well-posedness of the stochastic singular system in cases where the diffusion term is independent of the control variable. Moreover, the finite-horizon singular LQ Pareto optimal control problem was also studied. However, compared to its deterministic counterpart, the stochastic system presents new challenges in addition to the two crucial obstacles that arise from singular systems (further details are explained in the following paragraph). In our work, we aim to analyse the conditions under which the stochastic singular system, featuring both drift and diffusion terms dependent on control variables, has a unique solution. Furthermore, we investigate the singular LQ optimal control problem. To our best knowledge, our paper is the first one to deal with this kind of problem.
The main difficulties and our contributions are summarized as follows.
\begin{itemize}
\item Addressing the well-posedness of singular SDEs has remained a challenge. Except for \cite{gashi2013linear}, which explored the necessary conditions for regular backward systems, most relevant literature, including \cite{boukas2006stabilization, gao2013observer, ge2021impulse, huang2010stability, zhang2014stability, zhang2017linear, zhang2015some}, directly provided sufficient conditions without considering the necessary conditions. Inspired by \cite{dai1989singular, duan2010analysis}, which studied the deterministic case, we focus on the well-posedness of singular SDEs from a necessity perspective. Our approach is based on the Kronecker canonical form and the arbitrariness of control variables. Our results can be viewed as a generalization of \cite{dai1989singular, duan2010analysis} to the stochastic system. However, unlike the deterministic problem studied in  \cite{dai1989singular,
      duan2010analysis}, which only involves the decomposition of matrices $(E, A)$, our singular stochastic, where the state enters into the diffusion term,  involves the decomposition
      of three matrices $(E, A, C)$ due to the appearance of It\^o's integrals (see Section 3 of current paper for more details).  This feature makes that it is more challenging to establish the necessary conditions for the well-posedness of the singular SDE. By utilizing the properties of continuous martingales and constructing counterexamples, we tackle the complexities arising from stochastic integrals. Moreover, we argue that in \cite{dai1989singular,
      duan2010analysis}, the control variables were assumed to be sufficiently differentiable, which corresponds to It\^o processes in the stochastic counterpart. Unfortunately, the differentiability hypothesis is too strict. Instead, we propose the matrix symmetric technique to relax these control variable constraints. This technique is applicable not only in stochastic environments but also in deterministic cases, improving the results of \cite{dai1989singular, duan2010analysis}. Furthermore,
      our framework is more general than \cite{zhang2017linear}, underscoring the broader applicability of our conclusions. 
\item For the stochastic singular LQ control problem, it is impossible to directly establish the maximum principle because estimates on the trajectories cannot be obtained. To solve this problem, we use an equivalent transformation technique based on the Kronecker canonical form. This approach establishes a relationship between the singular LQ control problem and the normal LQ control problem. By solving the normal LQ control problem with the help of maximum principle, we derive the optimal control for the original problem. This method differs from the approach of completing the square in  \cite{zhang2014stability, zhang2017linear}. Moreover, we introduce the Riccati equation to give the feedback optimal control.
    \item In the study of the infinite-horizon normal LQ control problem, one need first to guarantee the finiteness of the cost functional, i.e. the non-emptiness of the admissible control set $\mathcal{U}{[0,\infty)}$ (for more details, see Section 5). It is usually utilize the definition of the stability to study the finiteness of the cost functional, see e.g. \cite{huang2015linear, rami2000linear, sun2018stochastic}. In contrast to these works,
         we employ the concept of stochastic exact controllability proposed by Peng \cite{peng1994backward} to ensure the finiteness of the cost functional. Additionally, we provide Popov-Belevitch-Hautus (PBH) rank criterion for assessing stochastic exact controllability, which differs from the algebraic criterion in \cite{peng1994backward}.
\item We give examples to illustrate the effectiveness of our results. It is noted that we provide numerical solutions for the relevant Riccati equations in Example \ref{example1}, while for Example \ref{example2}, we solve it explicitly.
       According to Theorem \ref{theorem5}, we derive a cubic equation that can be solved using the Cardano formula. Finally, by solving a linear SDE, we obtain the explicit optimal control, optimal state and the minimum value.
\end{itemize}

The present work is organized as follows. In Section 2, we list notations and preliminaries for the theory of matrices. In Section 3, we discuss the well-posedness of the stochastic singular system, which lays the foundation for the study of stochastic singular LQ control problem. In Section 4, we study the stochastic singular LQ control problem for finite-horizon case. In Section 5, we consider the stochastic singular LQ control problem for infinite-horizon case. In Section 6, we provide two examples. In section 7, we end the paper
with some concluding remarks. In the appendix, we complete the proof of Theorem \ref{theorem5}.
\section{Notations and preliminaries}
Throughout the paper, $\mathbb{C}$ represents the complex field. $\mathbb{R}^{n\times m}$ denotes the set of $n\times m$ real matrices. For any $A\in \mathbb{R}^{n\times m}$, $A^\mathrm{T}$ and $tr(A)$ denote respectively the transpose and the trace of $A$. We define $\vert A\vert=\sqrt{tr(A^\mathrm T A)}$ as the norm of matrix $A$. $A>(\geq)\ 0$ means that $A$ is positive definite (semi-definite) matrix. $I_n$ is the $n\times n$ unit matrix. $\mathbb{S}^n$ and $\mathbb{S}^n_+$ denote the set of $n\times n$ symmetric matrices and nonnegative definite matrices, respectively.

Suppose $(\Omega,\mathcal{F},\{\mathcal{F}_t\}_{t\geq 0},\mathbb{P})$ is a complete filtered probability space on which we define a one-dimensional standard Brownian motion $\{W(t), \ t\geq0\}$. We assume that $\{\mathcal{F}_t\}_{t\geq 0}$ is the natural filtration of $W$  augmented by all $\mathbb{P}$-null sets in $\mathcal{F}$.

Consider a finite time horizon $[0,T]$, where  $T > 0$ is fixed. We introduce the following notations.
\begin{itemize}
  \item $L^2_{\mathcal{F}_t}(0,T;\mathbb{R}^n)$: \ the space of $\mathcal{F}_t$-adapted $\mathbb{R}^n$-valued stochastic processes $x(\cdot)$ satisfying $\mathbb{E}\big[\int_0^T\vert x(t)\vert^2 dt\big] <\infty$;

  \item $L_{\mathcal{F}_T}^2(\Omega;\mathbb{R}^n)$: \ the space of all $\mathcal{F}_T$-measurable $\mathbb{R}^n$-valued random variable $\xi$ satisfying $\mathbb{E}\big[\vert\xi\vert^2 \big]< \infty$;
  \item
      $L_{\mathcal{F}_t}^2
      (\Omega;C([0,T];\mathbb{R}^n))$:
      \ the space of $\mathcal{F}_t$-adapted $\mathbb{R}^n$-valued continuous stochastic processes $x(\cdot)$ satisfying  $\mathbb{E}\big[\sup_{0\leq t\leq T}\vert x(t)\vert^2\big] <\infty$.
\end{itemize}

In the next step, we give some definitions and results for the matrix pencil and singular system, see \cite{dai1989singular,duan2010analysis,
gantmacher1974theory} for more details.
\begin{definition}{(Definition 1, Chapter XII, \cite{gantmacher1974theory})}
Given four matrices $G_1,\ E_1,\ G_2,\ E_2\in\mathbb{R}^{m\times n}$, a matrix pencil $G_1+\lambda E_1$ is called to be strictly equivalent to the other matrix pencil $G_2+\lambda E_2$ if there exist constant square invertible matrices $U\in\mathbb{R}^{m\times m}$ and $V\in\mathbb{R}^{n\times n}$ such that for any $\lambda\in\mathbb{C}$,
\begin{equation*}
U(G_1+\lambda E_1)V=G_2+\lambda E_2.
\end{equation*}
\end{definition}
According to the arbitrariness of $\lambda$, the above equality is equivalent to
$$UG_1V=G_2,\ UE_1V=E_2.
$$
\begin{definition}{(Definition 2, Chapter XII, \cite{gantmacher1974theory})}
\label{def3}
A matrix pencil $G_1+\lambda E_1$ is called regular if
\begin{itemize}
  \item[1)] $G_1$ and $E_1$ are square matrices of the same order $n$;
  \item[2)] There exists at least one $\lambda\in\mathbb{C}$ such that $\det(G_1+\lambda E_1)\neq0$.
\end{itemize}
\end{definition}
\begin{theorem}{(Lemma 1-2.2, \cite{dai1989singular}, or Theorem 3.1, \cite{duan2010analysis})}\label{Regular}
The square matrix pencil $G_1+\lambda E_1$ is regular if and only if there exist invertible matrices $U,V\in \mathbb{R}^{n\times n}$ such that
$$UE_1V=diag(I_{n_1},N),\  UG_1V=diag(\tilde{G}_1,I_{n_2}),
$$
where $n_1+n_2=n,\ \tilde{G}_1\in \mathbb{R}^{n_1\times n_1},\
N\in \mathbb{R}^{n_2\times n_2}$ and $N$ is a nilpotent matrix.
\end{theorem}
\begin{theorem}{(Theorem 5 (Kronecker canonical form), Chapter XII,  \cite{gantmacher1974theory})} \label{decop}
Every matrix pencil $G_1+\lambda E_1$ can be reduced (strictly equivalent) to a canonical quasi-diagonal form
\begin{equation}\label{kronecker}
\begin{aligned}
diag
\big(\mathbf{0}_{a\times b},
L_{\zeta_{b+1}},\cdots,
L_{\zeta_p},
L^\mathrm{T}_{\eta_{a+1}},
\cdots,L^\mathrm{T}_{\eta_q}, N^{(\rho_1)},\cdots,N^{(\rho_s)},
J+\lambda I\big),\\
\end{aligned}
\end{equation}
where
\begin{equation*}
\begin{aligned}
L_\zeta=\begin{pmatrix}
\lambda &1 &       &        &  \\
  &\lambda &1      &        &  \\
  &  &\ddots &\ddots  &  \\
  &  &       &\lambda       &1
  \end{pmatrix}
  \in\mathbb{R}^{\zeta\times
  (\zeta+1)}.
\end{aligned}
\end{equation*}

The choice of indices for $\zeta$ and $\eta$ is due to the fact that it is convenient here to take $\zeta_1=\zeta_2=
\cdots=\zeta_b=0$, $\eta_1=\eta_2=
\cdots=\eta_a=0$.
Besides,
$
N^{(\rho)}=I^{(\rho)}+\lambda H^{(\rho)}$. Here $I^{(\rho)}$ is a unit matrix of order $\rho$ and $H^{(\rho)}$ is a matrix of order $\rho$ whose elements in the first super-diagonal are 1, while the remaining elements are zero. $J$ is any matrix of normal form and $I$ is the unit matrix.
\end{theorem}
\begin{definition}{(Definition 3.1, \cite{duan2010analysis})}
\label{def2} If the matrix pencil $A+\lambda E$ in the system \eqref{ode} is regular, then we call the system \eqref{ode} regular.
\end{definition}
\begin{theorem}{(Theorem 2.7, \cite{duan2010analysis})}
The necessary and sufficient condition for a singular system \eqref{ode} to have a unique solution for any sufficiently differentiable $u(t)$ and some suitable initial values $x(0)\in\mathbb{R}^n$ is that the system is regular.
\end{theorem}

\section{The well-posedness of singular linear SDE}
In this section, we introduce a stochastic linear system, which is described by the following SDE
\begin{subnumcases}{\label{SDE}}
Edx(t) = [Ax(t)+Bu(t)]dt
+[Cx(t)+Du(t)]dW(t),
 \label{EQ_1}\\
Ex(0) = x_0,\label{EQ_2}
\end{subnumcases}
where $E,A,C\in \mathbb{R}^{n\times n}$, $B,D\in \mathbb{R}^{n\times r}$ are constant matrices; $x(\cdot)\in \mathbb{R}^n$ is the state vector, $u(\cdot)\in \mathbb{R}^r$ is the control process.

When $E$ is an invertible matrix, such a SDE \eqref{SDE} (also called system) is referred to as a normal SDE (stochastic normal system). In the following discussion, we are mainly concerned with the case where $E$ is singular, that is to say, rank $(E)=m<n$. Under the circumstance, SDE \eqref{SDE} (system) is called singular SDE (stochastic singular system). For convenience, we denote the stochastic singular system \eqref{SDE} by $(E,A,B,C,D)$. 

First of all, let us discuss some properties of stochastic singular systems.

\begin{remark}
The term ``singular'' is also commonly used in two other contexts. The first situation occurs when the system is driven by the following SDE
\begin{equation*}
\left\{
\begin{aligned}
dx(t)=&b(t,x(t),u(t))dt+\sigma(t,x(t),u(t))
dW(t)+\gamma(t)d\xi(t),\ t\in[0,T],\\
x(0)=&x_0,
\end{aligned}
\right.
\end{equation*}
where $\xi(\cdot)$ is a c\`{a}dl\`{a}g $\mathcal{F}_t$-adapted process which is of bounded variation and nondecreasing. Here $d\xi(\cdot)$ may be singular with respect to Lebesgue measure, thus $\xi(\cdot)$ is referred to as a singular control. The second situation arises when $\xi(\cdot)=0$ in the above equation, but $u(\cdot)$ trivially satisfies the first and second-order conditions of the maximum principle.\\
\indent In our paper, the term ``singular'' specifically refers to the appearance of a singular matrix on the left-hand side of the linear SDE. Therefore, the meanings of above three ``singular'' are different, it is simply a coincidental use of the term in different contexts within the control problem.
\end{remark}
\begin{definition}\label{def1}
If there exist another stochastic singular system $(\tilde{E},\tilde{A},\tilde{B},\tilde{C},
\tilde{D})$
\begin{equation}\label{SDE2}
\left\{\begin{aligned}
&\tilde{E}d\tilde{x}(t) = [\tilde{A}\tilde{x}(t)+\tilde{B}u(t)]dt
+[\tilde{C}\tilde{x}(t)+\tilde{D}u(t)]dW(t),\\
&\tilde{E}\tilde{x}(0) = x_0,
\end{aligned}
\right.
\end{equation}
and two invertible matrices $U,V\in \mathbb{R}^{n\times n}$ such that
$x=V\tilde{x},$
and $UEV= \tilde{E}, UAV=\tilde{A}, UCV=\tilde{C}, UB=\tilde{B}, UD=\tilde{D}$, then systems \eqref{SDE} and \eqref{SDE2} will be called restricted system equivalent (r.s.e).
\end{definition}
\begin{remark}
Definition \ref{def1} can be seen as an extension of the deterministic situation (see Definition 1 in  \cite{zhu2002singular}). Obviously, r.s.e is an equivalent relationship that possesses reflexivity, transitivity and invertibility.
\end{remark}

It is noted that the existence and uniqueness of equation \eqref{SDE} for all $u(\cdot)\in L^2_{\mathcal{F}_t}(0,T;\mathbb{R}^r)$ are significantly distinct from the normal SDE. The singularity of $E$ leads to the loss of information, so the well-posedness of the solution to equation \eqref{SDE} cannot be obtained directly by classical contraction mapping theorem. In the coming subsection, motivated by \cite{dai1989singular} and \cite{duan2010analysis}, we will give solvable conditions for equation \eqref{SDE} by the decomposition of Kronecker canonical form.
\subsection{The case of \ $C=\mathbf{0}$}
In this section, we suppose that $C=\mathbf{0}$. The singular SDE \eqref{SDE} now reduces to the subsequent simple form
\begin{equation}\label{SDE3}
\left\{
\begin{aligned}
&Edx(t) = [Ax(t)+Bu(t)]dt
+Du(t)dW(t),\\
&Ex(0) = x_0.
\end{aligned}
\right.
\end{equation}
According to Theorem \ref{decop}, we can derive the canonical quasi-diagonal form of $A+\lambda E$.
And then, we exchange the last two terms in \eqref{kronecker} by the row and column transformations. As a consequence, there exist two invertible matrices $M,N\in\mathbb{R}^{n\times n}$ such that
\begin{equation}\label{kronecker1}
\begin{aligned}
&\tilde{E}:= MEN=diag
(\mathbf{0}_{\tilde{n}_0\times n_0},
L_1,L_2,\ldots,L_{p},\check{L}_1,\check{L}_2,\ldots \check{L}_{q},I_h,G),\\
&\tilde{A}:= MAN=diag
(\mathbf{0}_{\tilde{n}_0\times n_0},
J_1,J_2,\ldots,J_{p},\check{J}_1,\check{J}_2,\ldots \check{J}_{q},A_1,I_g),
\end{aligned}
\end{equation}
where
\begin{equation*}
\mathbf{0}_{\tilde{n}_0\times n_0}\in\mathbb{R}
^{\tilde{n}_0\times n_0}, \quad
A_1\in \mathbb{R}^{h\times h},
\end{equation*}
\begin{equation}\label{matrix1}
\begin{aligned}
&L_i=\begin{pmatrix}
1 &0 &       &        &  \\
  &1 &0      &        &  \\
  &  &\ddots &\ddots  &  \\
  &  &       &1       &0
  \end{pmatrix},
&J_i=\begin{pmatrix}
0 &1 &       &        &  \\
  &0 &1      &        &  \\
  &  &\ddots &\ddots  &  \\
  &  &       &0       &1
  \end{pmatrix}
  \in\mathbb{R}^{\tilde{n}_i\times
  (\tilde{n}_i+1)},\\
\end{aligned}
\end{equation}
\begin{equation}\label{matrix2}
\begin{aligned}
&\check{L}_j=\begin{pmatrix}
1 &         &               &  \\
0 &1        &               &  \\
  &\ddots   &\ddots         &   \\
  &         &0              &1 \\
  &         &               &0
  \end{pmatrix},
&\check{J}_j=\begin{pmatrix}
0 &         &               &  \\
1 &0        &               &  \\
  &\ddots   &\ddots         &   \\
  &         &1              &0 \\
  &         &               &1
  \end{pmatrix}
  \in\mathbb{R}^{(n_j+1)\times
  n_j},
\end{aligned}
\end{equation}
\begin{equation*}
i=1,2,\cdots, p, \quad
j=1,2,\cdots,  q, \quad
\end{equation*}
\begin{equation}\label{matrix3}
\begin{aligned}
&\quad \ \ G=diag(G_1, G_2,\cdots, G_l)\in\mathbb{R}^{g\times g},\\
G_k=&\begin{pmatrix}
0 &1 &       &        &  \\
  &0 &1      &        &  \\
  &  &\ddots &\ddots  &  \\
  &  &       &0       &1  \\
  &  &       &        &0
  \end{pmatrix}
  \in\mathbb{R}^{m_k\times
  m_k},\
  k=1,2,\cdots,l.
\end{aligned}
\end{equation}
Thus the matrix $G$ is a nilpotent matrix. In fact, $G^{\tilde{m}}=\mathbf{0}_g$ with the index $\tilde{m}=\max\{m_k:k=1,2,\cdots,l\}$.
Moreover, the dimension of the above matrices has the following relationship
\begin{equation}\label{t21}
\left\{\begin{aligned}
&\tilde{n}_0+\sum_{i=1}^{p}\tilde{n}_i
+\sum_{j=1}^{q}(n_j+1)
+\sum_{k=1}^lm_k+h=n,\\
&n_0+\sum_{j=1}^{q}n_j
+\sum_{i=1}^{p}(\tilde{n}_i+1)
+\sum_{k=1}^lm_k+h=n,\\
&\sum_{k=1}^lm_k=g.
\end{aligned}
\right.
\end{equation}

For avoiding heavy notations, we can omit the time variable $t$ if no confusion arises.
Let
$x=N\tilde{x},
$
and multiplying \eqref{SDE3} by the invertible matrix $M$ on the left, we have
\begin{equation}\label{SDE4}
\left\{\begin{aligned}
&\tilde{E}d\tilde{x}(t) = [\tilde{A}\tilde{x}(t)+\tilde{B}u(t)]dt
+\tilde{D}u(t)dW(t),\\
&\tilde{E}\tilde{x}(0) = Mx_0,
\end{aligned}
\right.
\end{equation}
where $\tilde{E}$ and $\tilde{A}$ are given by \eqref{kronecker1}. According to the structure of $\tilde{E}$ and $\tilde{A}$, the matrices $\tilde{B}$ and $\tilde{D}$ can be represented as the following form, respectively,
\begin{equation}\label{T1}
\tilde{B}=MB=
\left[
\begin{array}{c;{3pt/1pt}c;{3pt/1pt}c
;{3pt/1pt}c
;{3pt/1pt}c}
B_{1}^\mathrm{T}
& \begin{matrix}
B_{21}^\mathrm{T}
&\cdots &B_{2p}^\mathrm{T}
\end{matrix}&
\begin{matrix}
B_{31}^\mathrm{T}
&\cdots
&B_{3q}^\mathrm{T}
\end{matrix}
&B_4^\mathrm{T}
&B_5^\mathrm{T}
\end{array}
\right]^\mathrm{T},
\end{equation}

\begin{equation}\label{T2}
\tilde{D}=MD=
\left[
\begin{array}{c;{3pt/1pt}c;{3pt/1pt}c
;{3pt/1pt}c
;{3pt/1pt}c}
D_{1}^\mathrm{T}
& \begin{matrix}
D_{21}^\mathrm{T}
&\cdots &D_{2p}^\mathrm{T}
\end{matrix}&
\begin{matrix}
D_{31}^\mathrm{T}
&\cdots
&D_{3q}^\mathrm{T}
\end{matrix}
&D_4^\mathrm{T}
&D_5^\mathrm{T}
\end{array}
\right]^\mathrm{T},
\end{equation}
Meanwhile, the state vector $\tilde{x}$ has the similar decomposition
\begin{equation}\label{t13}
\tilde{x}=
\left[
\begin{array}{c;{3pt/1pt}c;{3pt/1pt}c
;{3pt/1pt}c
;{3pt/1pt}c}
\tilde{x}_1^\mathrm{T}
& \begin{matrix}
\tilde{x}_{21}^\mathrm{T} 
&\cdots &\tilde{x}_{2p}^\mathrm{T}
\end{matrix}&
\begin{matrix}
\tilde{x}_{31}^\mathrm{T}
&\cdots
&\tilde{x}_{3q}^\mathrm{T}
\end{matrix}
&\tilde{x}_4^\mathrm{T}
&\tilde{x}_5^\mathrm{T}
\end{array}
\right]^\mathrm{T}.
\end{equation}
According to transformations \eqref{kronecker1} and \eqref{T1}-\eqref{t13}, equation \eqref{SDE4} can be rewritten as
\begin{subequations}\label{AC0}
 \begin{align}
  &\mathbf{0}_{\tilde{n}_0\times n_0}d\tilde{x}_{1}(t)
  =[\mathbf{0}_{\tilde{n}_0\times n_0}\tilde{x}_{1}(t)
  +B_{1}u(t)]dt
  +D_{1}u(t)dW(t), \label{SDE5a}\\
  &L_id\tilde{x}_{2i}(t)
  =[J_i\tilde{x}_{2i}(t)+B_{2i}u(t)]dt
  +D_{2i}u(t)dW(t),
  \quad i=1,2,\cdots,p ,\label{SDE5b}\\
  &\check{L}_jd\tilde{x}_{3j}(t)
  =[\check{J}_j\tilde{x}_{3j}(t)
  +B_{3j}u(t)]dt
   +D_{3j}u(t)dW(t),
  \quad j=1,2,\cdots,q, \label{SDE5c}\\
  &d\tilde{x}_4(t)
  =[A_1\tilde{x}_4(t)+B_4u(t)]dt
   +D_4u(t)dW(t), \label{SDE5d}\\
  &G_kd\tilde{x}_{5k}(t)
  =[\tilde{x}_{5k}(t)
  +B_{5k}u(t)]dt
  +D_{5k}u(t)dW(t), \quad k=1,2,\cdots,l, \label{SDE5e}
  \end{align}
\end{subequations}
where $B_5,\ D_5$ and $\tilde{x}_5$ have the following decomposition, respectively
\begin{equation*}
\begin{aligned}
&B_5=(B_{51}^\mathrm{T},
B_{52}^\mathrm{T},\cdots,
B_{5l}^\mathrm{T})^\mathrm{T},\
D_5=(D_{51}^\mathrm{T},
D_{52}^\mathrm{T},\cdots,
D_{5l}^\mathrm{T})^\mathrm{T},\
\tilde{x}_5=\big(
\tilde{x}_{51}^\mathrm{T},
\tilde{x}_{52}^\mathrm{T},\cdots,
\tilde{x}_{5l}^\mathrm{T}\big)
^\mathrm{T},\\
&B_{5k},\ D_{5k} \in\mathbb{R}^{m_k\times r}, \ \tilde{x}_{5k}(\cdot)\in\mathbb{R}^{m_k},
\ 1\leq k\leq l.
\end{aligned}
\end{equation*}
Therefore, the solution to equation \eqref{SDE3} or the solution to equation \eqref{SDE4} is determined by equations \eqref{SDE5a}-\eqref{SDE5e} and the initial condition. Now, let us verify the solvability of each equation.
\begin{enumerate}
  \item Solvability of \eqref{SDE5a}. \ We have
      $\mathbf{0}=B_1u(t)dt+D_1u(t)dW(t).
      $
      By the property of continuous martingale (see Proposition 1.2 in Chapter IV, \cite{revuz2013continuous}), we have $B_1u(t)=\mathbf{0},\ D_1u(t)=\mathbf{0}$, for a.e. $t$, then equation \eqref{SDE5a} can be solved for any control $u(\cdot)\in L^2_{\mathcal{F}_t}
      (0,T;\mathbb{R}^r)$ only when $B_{1}=D_{1}=\mathbf{0}_{\tilde{n}_0\times r}$. In this case, equation \eqref{SDE5a} is an identical equation. Therefore, $\tilde{x}_{1}(\cdot)$ has either infinite number of solutions or no solution.
  \item  Solvability of \eqref{SDE5b}. We define state vector $\tilde{x}_{2i}^\mathrm{T}
      =(\tilde{x}_{2i,1},
      \cdots,\tilde{x}_{2i,\tilde{n}_i},
\tilde{x}_{2i,\tilde{n}_i+1})$ and denote $B_{2i}^\mathrm{T}=\big(
      B_{2i,1}^\mathrm{T},
      \cdots,B_{2i,\tilde{n}_i}^\mathrm{T}
      \big)$ and $D_{2i}^\mathrm{T}=\big(
      D_{2i,1}^\mathrm{T},
      \cdots,
      D_{2i,\tilde{n}_i}^\mathrm{T}
      \big)$, for $i=1,2,\cdots,p$. Equation \eqref{SDE5b} is constituted by a collection of equations of the following form
\begin{equation*}
\begin{aligned}
d
\begin{pmatrix}
1 &0 &       &        &  \\
  &1 &0      &        &  \\
  &  &\ddots &\ddots  &  \\
  &  &       &1       &0
  \end{pmatrix}\tilde{x}_{2i}(t)=&\left[ \begin{pmatrix}
0 &1 &       &        &  \\
  &0 &1      &        &  \\
  &  &\ddots &\ddots  &  \\
  &  &       &0       &1
  \end{pmatrix}\tilde{x}_{2i}(t)+
  \begin{pmatrix}
  B_{2i,1}\\
  B_{2i,2}\\
  \vdots\\
  B_{2i,\tilde{n}_i}
  \end{pmatrix}u(t) \right] dt\\
  &+\begin{pmatrix}
  D_{2i,1}\\
  D_{2i,2}\\
  \vdots\\
  D_{2i,\tilde{n}_i}
  \end{pmatrix}u(t)dW(t),\ i=1,2,\cdots,p.
  \end{aligned}
\end{equation*}
It is obvious that the above equation is equivalent to
\begin{equation*}
\left\{\begin{aligned}
d\tilde{x}_{2i,1}(t)=&[
\tilde{x}_{2i,2}(t)+B_{2i,1}u(t)]dt
+D_{2i,1}u(t)dW(t),\\
d\tilde{x}_{2i,2}(t)=&[
\tilde{x}_{2i,3}(t)+B_{2i,2}u(t)]dt
+D_{2i,2}u(t)dW(t),\\
&\vdots\\
d\tilde{x}_{2i,\tilde{n}_i}(t)
=&[\tilde{x}_{2i,\tilde{n}_i+1}(t)
+B_{2i,\tilde{n}_i}u(t)]dt
+D_{2i,\tilde{n}_i}u(t)dW(t).\\
\end{aligned}
\right.
\end{equation*}
When $\tilde{x}_{2i,\tilde{n}_i+1}$ and the initial condition are given, $\tilde{x}_{2i,\tilde{n}_i}$ has a unique solution. And then, $\tilde{x}_{2i,\tilde{n}_i-1}, \cdots, \tilde{x}_{2i,2}, \tilde{x}_{2i,1}$ can be solved uniquely. That is to say, \eqref{SDE5b} has one free unknown $\tilde{x}_{2i,\tilde{n}_i+1}$. Thus, the solution to \eqref{SDE5b} is not unique.
  \item Solvability of \eqref{SDE5c}. \ Let $\tilde{x}_{3j}
      ^\mathrm{T}=
      (\tilde{x}_{3j,1},
      \cdots,
\tilde{x}_{3j,n_j})$,  $B_{3j}^\mathrm{T}=\big(
      B_{3j,1}^\mathrm{T},
      \cdots,B_{3j,n_j}^\mathrm{T},
      B_{3j,n_j+1}^\mathrm{T}\big)$ and\\
      $D_{3j}^\mathrm{T}=\big(
      D_{3j,1}^\mathrm{T},\cdots, D_{3j,n_j}^\mathrm{T},
      D_{3j,n_j+1}^\mathrm{T}\big)$, for $j=1,2,\cdots,q$. Every equation in \eqref{SDE5c} is of the following form
\begin{equation*}
\begin{aligned}
d\begin{pmatrix}
1 &         &               &  \\
0 &1        &               &  \\
  &\ddots   &\ddots         &   \\
  &         & 0             &1 \\
  &         &               &0
  \end{pmatrix}\tilde{x}_{3j}(t)=&\left[
  \begin{pmatrix}
0 &         &               &  \\
1 &0        &               &  \\
  &\ddots   &\ddots         &   \\
  &         &1              &0 \\
  &         &               &1
  \end{pmatrix}\tilde{x}_{3j}(t)+
  \begin{pmatrix}
  B_{3j,1}\\
  B_{3j,2}\\
  \vdots\\
  B_{3j,n_j}\\
  B_{3j,n_j+1}
  \end{pmatrix}u(t)\right] dt\\
  &+\begin{pmatrix}
  D_{3j,1}\\
  D_{3j,2}\\
  \vdots\\
  D_{3j,n_j}\\
  D_{3j,n_j+1}
  \end{pmatrix}u(t)dW(t),\ j=1,2,\cdots,q.
\end{aligned}
\end{equation*}

Then, we get
\begin{subnumcases}{\label{SDE19}}
d\tilde{x}_{3j,1}(t)
=B_{3j,1}u(t)dt+D_{3j,1}u(t)dW(t), \label{SDE19-1}\\
d\tilde{x}_{3j,2}(t)
=[\tilde{x}_{3j,1}(t)
+B_{3j,2}u(t)]dt
+D_{3j,2}u(t)dW(t),
\label{SDE19-2}\\
\qquad \qquad \vdots \notag\\
d\tilde{x}_{3j,n_j}(t)
=[\tilde{x}_{3j,n_j-1}(t)+B_{3j, n_j}u(t)]dt
+D_{3j,n_j}u(t)dW(t),
\label{SDE19-3}\\
0=[\tilde{x}_{3j,n_j}(t)+B_{3j, n_j+1}u(t)]dt
+D_{3j,n_j+1}u(t)dW(t).
\label{SDE19-4}
\end{subnumcases}
  On the one hand, noticing that $\tilde{x}_{3j,1}$ can be solved uniquely in equation \eqref{SDE19-1} for each fixed $u(\cdot)\in L^2_{\mathcal{F}_t}
      (0,T;\mathbb{R}^r)$. Substituting $\tilde{x}_{3j,1}$ into equation \eqref{SDE19-2}, we derive a unique $\tilde{x}_{3j,2}$. Iterating the above procedures, we solve $\tilde{x}_{3j,1}, \tilde{x}_{3j,2},\cdots, \tilde{x}_{3j,n_j}$ uniquely. On the other hand,
  according to equation \eqref{SDE19-4} and the property of continuous martingale, we obtain for a.e. $t$,
\begin{equation}\label{t22}
\left\{\begin{aligned}
&\tilde{x}_{3j,n_j}(t)+B_{3j,n_j+1}u(t)=0, \\
&D_{3j,n_j+1}u(t)=0.
\end{aligned}
\right.
\end{equation}
Since equality \eqref{t22} holds for any control $u(\cdot)\in L^2_{\mathcal{F}_t}(0,T;\mathbb{R}^r)$, we derive $D_{3j,n_j+1}=\mathbf{0}_{r}$.
Denote by $\Phi_{1j,n_{j}+1}=B_{3j,n_j+1}
^\mathrm{T}B_{3j,n_j+1}$.
For any $B_{3j,n_j+1}^\mathrm{T}\in\mathbb{R}^r$, rank $(\Phi_{1j,n_{j}+1})\leq1$ because $B_{3j,n_j+1}^\mathrm{T}$ is a column vector. In the next step, we discuss two cases regarding to the above equation. And we call this method the matrix symmetric technique.
\begin{itemize}
  \item[1)] $B_{3j,n_j+1}^\mathrm{T}\neq \mathbf{0}_{r\times1}$. Namely, there exist invertible matrices $P_{1j,n_{j}+1}$ and $Q_{1j,n_{j}+1}\in\mathbb{R}
      ^{r\times r}$ such that
      $$\Phi_{1j,n_{j}+1}=
      P_{1j,n_{j}+1}\begin{pmatrix}
      1   &0 &\cdots &0  \\
      0  &0 &\cdots &0 \\
      \vdots &\vdots   & &\vdots\\
      0 &0 &\cdots &0
      \end{pmatrix}Q_{1j,n_{j}+1}.$$
      Substituting $\tilde{x}_{3j,n_j}
      =-B_{3j,n_j+1}u$ into equation \eqref{SDE19-3} and let $\tilde{u}=Q_{1j,n_{j}+1}u$, we get
      \begin{equation}\label{u2}
      \begin{aligned}
      &d\begin{pmatrix}
      1   &0 &\cdots &0  \\
      0  &0 &\cdots &0 \\
      \vdots &\vdots   & &\vdots\\
      0 &0 &\cdots &0
      \end{pmatrix}\tilde{u}(t)
      =-P_{1j,n_{j}+1}^{-1}B_{3j,
      n_j+1}^\mathrm{T}\big[
      \tilde{x}_{3j,n_j-1}(t)\\
      &+B_{3j,n_j}Q_{1j,n_{j}+1}^{-1}\tilde{u}(t)\big]dt
      -P_{1j,n_{j}+1}^{-1}B^\mathrm{T}_{
      3j,n_j+1}D_{3j,
      n_j}Q_{1j,n_{j}+1}^{-1}\tilde{u}(t)dW(t).
      \end{aligned}
      \end{equation}
      It is easy to verify that equation \eqref{u2} does not hold for $\tilde{u}(t)=(W(t),\cdots,
      W(t))^\mathrm{T}$. Thus, equation \eqref{u2} does not have a unique solution for any
      $u(\cdot)\in L^2_{\mathcal{F}_t}(0,T;
      \mathbb{R}^r)$.
  \item[2)] $B_{3j,n_j+1}^\mathrm{T}
      =\mathbf{0}_{r\times1}$.
      From \eqref{t22} and \eqref{SDE19-3}, we conclude that $\tilde{x}_{3j,nj}(t)\equiv0$. Similarly, by  \eqref{SDE19-3}, we get for a.e. $t$
      \begin{equation*}
\left\{\begin{aligned}
&\tilde{x}_{3j,n_j-1}(t)+B_{3j, n_j}u(t)=0, \\
&D_{3j,n_j}u(t)=0.
\end{aligned}
\right.
\end{equation*}
\end{itemize}

Iterating the above procedures, equation \eqref{SDE5c} has a unique solution $\tilde{x}_{3j}$ for any $u(\cdot)\in L^2_{\mathcal{F}_t}(0,T;
      \mathbb{R}^r)$ only when $B_{3j}=D_{3j}=\mathbf{0}_{(n_j+1)\times r}$ and the initial value $\tilde{x}_{3j}(0)=
      \mathbf{0}_{n_j\times1}$, for $j=1,2,\cdots,q$. Moreover, the unique solution is $\tilde{x}_{3j}
      \equiv\mathbf{0}_{n_j\times1},\ j=1,2,\cdots,q$.
  \item Solvability of \eqref{SDE5d}. Equation \eqref{SDE5d} is a linear SDE, which has a unique solution for any $u(\cdot)\in L^2_{\mathcal{F}_t}
      (0,T;\mathbb{R}^r)$ and any initial value $\tilde{x}_4(0)$.
  \item  Solvability of \eqref{SDE5e}. Denote $\tilde{x}_{5k}
      ^\mathrm{T}=
      (\tilde{x}_{5k,1},
      \tilde{x}_{5k,2},\cdots,
\tilde{x}_{5k,m_k})$, equation \eqref{SDE5e} is composed of the coming equation
\begin{equation*}
d\begin{pmatrix}
0 &1 &       &        &  \\
  &0 &1      &        &  \\
  &  &\ddots &\ddots  &  \\
  &  &       &0       &1  \\
  &  &       &        &0
  \end{pmatrix}\tilde{x}_{5k}(t)
  =\left[\tilde{x}_{5k}(t)+
 \begin{pmatrix}
  B_{5k,1}\\
  B_{5k,2}\\
  \vdots\\
  B_{5k,m_k-1}\\
  B_{5k,m_k}
  \end{pmatrix}u(t)\right]dt
  +\begin{pmatrix}
  D_{5k,1}\\
  D_{5k,2}\\
  \vdots\\
  D_{5k,m_k-1}\\
  D_{5k,m_k}
  \end{pmatrix}u(t)dW(t)
\end{equation*}
where $B_{5k}^\mathrm{T}=
      (B_{5k,1}^\mathrm{T},
      B_{5k,2}^\mathrm{T},
      \cdots,B_{5k,m_k}^\mathrm{T})
      $,
      $D_{5k}^\mathrm{T}=
      (D_{5k,1}^\mathrm{T},
      D_{5k,2}^\mathrm{T},
      \cdots,D_{5k,m_k}^\mathrm{T})
      $. Then, we obtain
\begin{subnumcases}{\label{SDE14}}
d\tilde{x}_{5k,2}(t)=[
\tilde{x}_{5k,1}(t)+B_{5k,1}u(t)]dt
+D_{5k,1}u(t)dW(t),
\label{EQ14_1}\\
d\tilde{x}_{5k,3}(t)=[
\tilde{x}_{5k,2}(t)+B_{5k,2}u(t)]dt
+D_{5k,2}u(t)dW(t),
\label{EQ14_2}\\
\qquad \qquad \vdots \notag\\
d\tilde{x}_{5k,m_k}(t)=
[\tilde{x}_{5k,m_k-1}(t)+B_{ 5k,m_k-1}u(t)]dt
+D_{5k,m_k-1}u(t)dW(t), \label{EQ14_3}\\
0=[\tilde{x}_{5k,m_k}(t)+B_{5k,m_k}u(t)]dt
+D_{5k,m_k}u(t)dW(t). \label{EQ14_4}
\end{subnumcases}
From equation \eqref{EQ14_4}, we obtain for a.e. $t$,
\begin{equation}\label{t1}
\left\{\begin{aligned}
&\tilde{x}_{5k,m_k}(t)=-B_{5k,m_k}u(t), \\
&D_{5k,m_k}=\mathbf{0}_{1\times r}.
\end{aligned}
\right.
\end{equation}
Similarly, we discuss the solvability of \eqref{t1} by the matrix symmetric technique. Define $\Phi_{2k,m_k}
=B^\mathrm{T}_{5k,m_k}B_{5k,m_k}$.
\begin{itemize}
  \item[1)] $B_{5k,m_k}^\mathrm{T}
      \neq\mathbf{0}_{r\times1}$. Namely, there exist invertible matrices $P_{2k,m_k}$ and $Q_{2k,m_k}\in \mathbb{R}^{r\times r}$ such that
      $$\Phi_{2k,m_k} =P_{2k,m_k}
      \begin{pmatrix}
      1   &0 &\cdots &0  \\
      0  &0 &\cdots &0 \\
      \vdots &\vdots   & &\vdots\\
      0 &0 &\cdots &0
      \end{pmatrix}Q_{2k,m_k}.$$
      Let $\tilde{u}=Q_{2k,m_k}u$, we get
      \begin{equation}\label{u3}
      \begin{aligned}
      &d\begin{pmatrix}
      1   &0 &\cdots &0  \\
      0  &0 &\cdots &0 \\
      \vdots &\vdots   & &\vdots\\
      0 &0 &\cdots &0
      \end{pmatrix}\tilde{u}(t)
      =-P_{2k,m_k}^{-1}B^\mathrm{T}
      _{5k,m_k}\big[
      \tilde{x}_{5k,m_k-1}(t)\\
      &+B_{5k,m_k-1}Q_{2k,m_k}^{-1}
      \tilde{u}(t)\big]dt
      -P_{2k,m_k}^{-1}
      B^\mathrm{T}_ {5k,m_k}D_{5k,m_k-1}
      Q_{2k,m_k}^{-1}\tilde{u}(t)
      dW(t).
      \end{aligned}
      \end{equation}
      Similar to the discussion of \eqref{u2}, equation \eqref{u3} can not have a solution for any $u(\cdot)\in   L^2_{\mathcal{F}_t}
      (0,T;\mathbb{R}^r)$.
  \item[2)] $B_{5k,m_k}^\mathrm{T}
      =\mathbf{0}_{r\times1}
      $. Then, from \eqref{EQ14_3} and \eqref{t1}, we have  $\tilde{x}_{5k,m_k}(t)\equiv0$ and for a.e. $t$
      \begin{equation*}
\left\{\begin{aligned}
\tilde{x}_{5k,m_k-1}(t)&=-B_{5k, m_k-1}u(t),\\
D_{5k,m_k-1}&=\mathbf{0}_{1\times r}.
\end{aligned}
\right.
\end{equation*}
\end{itemize}
Iterating the above procedures, equation \eqref{SDE5e} has a unique solution for any $u(\cdot)\in L^2_{\mathcal{F}_t}(0,T;\mathbb{R}^r)$ only when $B_{5k,2}=\cdots=B_{ 5k,m_k-1}=B_{5k,m_k}=\mathbf{0}_{
1\times r}$, $D_{5k}=\mathbf{0}_{
m_k\times r}$ and the initial values $\tilde{x}_{5k,2}(0)
=\cdots=\tilde{x}_{5k,m_k-1}(0)
=\tilde{x}_{5k,m_k}(0)=0$, $k=1,2,\cdots,l$. Moreover, the unique solution to equation \eqref{SDE5e} is $\left(-B_{5k,1}u,0,\cdots,0\right)
^\mathrm{T}$, $k=1,2,\cdots,l$.
\end{enumerate}

It is noted that $E$ and $A$ are $n\times n$ square matrices. According to relationship \eqref{t21}, if equations \eqref{SDE5a} and \eqref{SDE5b} vanish, equation \eqref{SDE5c} must not exist too.
Therefore, in order to ensure the existence and uniqueness of the solution to equation \eqref{SDE3} for any $u(\cdot)\in L^2_{\mathcal{F}_t}(0,T;\mathbb{R}^r)$, \eqref{SDE5a}-\eqref{SDE5c} should disappear and the following conditions hold, for $k=1,2,\cdots,l$,
\begin{equation*}
\begin{aligned}
&B_{5k,2}=\cdots=B_{5k,m_k-1}=B_{5k,m_k}
=\mathbf{0}_{1\times r},\ D_{5k}=\mathbf{0}_{m_k\times r},\notag\\
&\tilde{x}_{5k,2}(0)
=\cdots=\tilde{x}_{5k,m_k-1}(0)
=\tilde{x}_{5k,m_k}(0)=0.
\end{aligned}
\end{equation*}
Noting that the above initial condition is equivalent to $G\tilde{x}_{5}(0)
=\mathbf{0}_{g\times1}$. Then, from \eqref{SDE4}, we obtain that $\begin{pmatrix}
\mathbf{0}_{(n-h)\times h}
&I_{n-h}
\end{pmatrix}Mx_0
=\mathbf{0}_{(n-h)\times1}$ is a necessary initial condition for the well-posedness of equation \eqref{SDE3}. Therefore, based on the above discussion, we get the following theorem which extends the results of deterministic singular system (see Theorem 2.7 in \cite{duan2010analysis} or Theorem 2.3 in current paper).
\begin{theorem}\label{C0}
Equation \eqref{SDE3} has a unique solution $x(\cdot)\in L^2_{\mathcal{F}_t}(0,T;\mathbb{R}^n)$ for any $u(\cdot)\in L^2_{\mathcal{F}_t}(0,T;\mathbb{R}^r)$ if and only if there exist invertible matrices $M$ and $N\in \mathbb{R}^{n\times n}$ such that
\begin{equation*}\label{t19}
\begin{aligned}
&MEN=diag(I_{h},G), \quad MAN=diag(A_1,I_{n-h}),\\
&MB=\begin{pmatrix}
B_1\\
B_2
\end{pmatrix},\quad
MD=\begin{pmatrix}
D_1\\
\mathbf{0}_{(n-h)\times r}
\end{pmatrix}.
\end{aligned}
\end{equation*}
where
$B_2=\big(B_{21}^\mathrm{T},
B_{22}^\mathrm{T},
\cdots,
B_{2l}^\mathrm{T}\big)^\mathrm{T}$
and
$$B_{2k}=
\begin{pmatrix}
b_{k1}\\
0\\
\vdots\\
0
\end{pmatrix}\in\mathbb{R}^{m_k\times r},\ 1\leq k \leq l,$$
$G \in \mathbb{R}^{(n-h)\times (n-h)}$ is a nilpotent matrix of form \eqref{matrix3},
$A_1\in\mathbb{R}^{h\times h}$ is any matrix,
$B_1,D_1\in\mathbb{R}^{h\times r},
B_2\in\mathbb{R}^{(n-h)\times r}$ and $b_{k1}\in \mathbb{R}^{1\times r},\
1\leq k\leq l$. In addition,
$\sum_{k=1}^lm_k=n-h$.
Moreover, the initial value satisfies
$\begin{pmatrix}
\mathbf{0}_{(n-h)\times h}
&I_{n-h}
\end{pmatrix}Mx_0
=\mathbf{0}_{(n-h)\times1}.$
\end{theorem}
\noindent{\bf Proof}\quad
To summarize the above analysis, we derive the necessity, and we now prove the sufficiency.
Let $x=N\begin{pmatrix}
x_1\\
x_2
\end{pmatrix}$, $x_2^\mathrm{T}=(x_{21}^\mathrm{T},\cdots x_{2l}^\mathrm{T})^\mathrm{T}$ and then we obtain that equation \eqref{SDE3} is restricted equivalent to the following
\begin{equation*}
\left\{
\begin{aligned}
dx_1(t)&=[A_1x_1(t)+B_1u(t)]dt
+D_1u(t)dW(t),\\
dG_kx_{2k}(t)&=
[x_{2k}(t)+B_{2k}u(t)]dt,\\
x_1(0)&=\begin{pmatrix}
I_h &\mathbf{0}_{h\times(n-h)}
\end{pmatrix}Mx_0,\ G_kx_{2k}(0)=\mathbf{0}_{m_k},\ k=1,2,\cdots,l.
\end{aligned}
\right.
\end{equation*}
Let $x_{2k}^\mathrm{T}=(x_{2k,1},\cdots
x_{2k,m_k})$. From \eqref{matrix3}, we have
\begin{equation*}
\left\{
\begin{aligned}
&dx_1(t)=[A_1x_1(t)+B_1u(t)]dt
+D_1u(t)dW(t),\\ &x_1(0)=\begin{pmatrix}
I_h &\mathbf{0}_{h\times(n-h)}
\end{pmatrix}Mx_0,
\end{aligned}
\right.
\end{equation*}
and for $k=1,2,\cdots,l$,
\begin{equation*}
\left\{
\begin{aligned}
&dx_{2k,2}(t)=
[x_{2k,1}(t)+b_{k1}u(t)]dt,\\
&\qquad \qquad \vdots\\
&dx_{2k,m_k}(t)=x_{2k,m_k-1}(t)dt\\
&0=x_{2k,m_k}(t)dt,\\
&x_{2k,2}(0)=\cdots
=x_{2k,m_k}(0)=0,
\end{aligned}
\right.
\end{equation*}
which implies that euqation \eqref{SDE3} has a unique solution for any $u(\cdot)\in L^2_{\mathcal{F}_t}(0,T;\mathbb{R}^r)$. This completes the proof.
\hfill$\Box$
\begin{remark}
It can be seen that the matrix symmetric technique help us remove the assumption about differentiability of control in \cite{dai1989singular} and \cite{duan2010analysis}. The necessary and sufficient condition in this paper looks stronger than the one in \cite{dai1989singular} and \cite{duan2010analysis}. This is caused by the fact that the control set we consider here is $L^2_{\mathcal{F}_t}(0,T;
\mathbb{R}^r)$, and it is larger than the set of It\^o processes which are reduced to differential functions in deterministic case.
\end{remark}
\subsection{The case of \ $C\neq\mathbf{0}$}
When $C\neq\mathbf{0}$, obtaining the necessary and sufficient conditions for the existence and uniqueness of equation \eqref{SDE} is challenging because it appears necessary to simultaneously transform the three matrices $E,A,C$ if using similar arguments as in Section 3.1.
To solve this problem, motivated by Theorem \ref{C0}, we will give the definition of strongly regular in the following and we will discuss the well-posedness of equation \eqref{SDE} for any $u(\cdot)\in L^2_{\mathcal{F}_t}(0,T;\mathbb{R}^r)$ under strongly regular condition. First, similar to Definition \ref{def2}, we call the stochastic system $(E,A,B,C,D)$ regular if the matrix pencil $A+\lambda E$ is regular, see \cite{gao2013observer, huang2010stability, zhang2014stability}.
\begin{definition}\label{SRegular}
A controlled system $(E,A,B,C,D)$ is called strongly regular if there exist invertible matrices $M$ and $N$ $\in \mathbb{R}^{n\times n}$ such that
\begin{itemize}
  \item[(1)] $MEN=diag(I_h,G),\
      MAN=diag(A_1,I_{n-h}),$
  \item[(2)] $MB=\begin{pmatrix}
B_1\\
B_2
\end{pmatrix},$
\end{itemize}
where
$G \in \mathbb{R}^{(n-h)\times (n-h)}$ is a nilpotent matrix of form \eqref{matrix3},
$A_1\in\mathbb{R}^{h\times h}$ is any matrix,
$B_1\in\mathbb{R}^{h\times r},
B_2\in\mathbb{R}^{(n-h)\times r}$,
$B_2=\big(B_{21}^\mathrm{T},
B_{22}^\mathrm{T},
\cdots,
B_{2l}^\mathrm{T}\big)^\mathrm{T}$
and
$$B_{2k}=
\begin{pmatrix}
b_{k1}\\
0\\
\vdots\\
0
\end{pmatrix}\in\mathbb{R}^{m_k\times r},\ b_{k1}\in \mathbb{R}^{1\times r}, \ 1\leq k \leq l.$$
In addition, the dimensions of the above matrices have the relationship
$\sum_{k=1}^lm_k=n-h$.
\end{definition}


The system $(E,A,B,C,D)$ being regular is equivalent to the condition stated in the first term of Definition \ref{SRegular} (see Theorem 2.1). It can be seen that the definition of strongly regular includes also constraints on the matrix $B$. Noting that in the case of $C=\mathbf{0}$, strongly regular is a necessary condition for the existence and  uniqueness of the solution to the system. However, when $C\neq\mathbf{0}$, the situation becomes different. Indeed, even if the system is not regular, the system may still have a unique solution. In another perspective, the state variable $x(\cdot)$ can be viewed as a ``good'' factor when it enters into the diffusion term of the system, and the condition for the well-posedness of singular SDEs can be weaken.
\begin{example}\label{ex1}
Consider the controlled singular SDE
\begin{equation*}
\left\{
\begin{aligned}
&d\begin{pmatrix}
0 &0 &0\\
0 &0 &1\\
0 &0 &0
\end{pmatrix}
\begin{pmatrix}
x_1\\
x_2\\
x_3
\end{pmatrix}
=
\left[
\begin{pmatrix}
0 &0 &0\\
0 &1 &0\\
0 &0 &1
\end{pmatrix}
\begin{pmatrix}
x_1\\
x_2\\
x_3
\end{pmatrix}
+\begin{pmatrix}
0\\
1\\
0
\end{pmatrix}
u \right]
dt+
\left[
\begin{pmatrix}
1 &0 &0\\
0 &1 &0\\
0 &0 &1
\end{pmatrix}
\begin{pmatrix}
x_1\\
x_2\\
x_3
\end{pmatrix}
+\begin{pmatrix}
1\\
1\\
0
\end{pmatrix}
u\right]dW(t),\\
&x_1(0)=0,\ x_2(0)=0,\ x_3(0)=0.
\end{aligned}
\right.
\end{equation*}
which is equivalent to
\begin{equation*}
\left\{\begin{aligned}
&0=(x_1+u)dW(t),\\
&dx_3=(x_2+u)dt+(x_2+u)dW(t),\\
&0=x_3dt+x_3dW(t),\\
&x_1(0)=0,\ x_2(0)=0,\ x_3(0)=0.
\end{aligned}
\right.
\end{equation*}

By solving the above equations, the system has a unique solution $x_1=x_2=-u$, a.e., $x_3\equiv0$. It can be verified that the above singular SDE is not regular (see Definition \ref{def3} and Definition \ref{def2}). Indeed, $\det\left\vert
\begin{pmatrix}
0 &0 &0\\
0 &1 &0\\
0 &0 &1
\end{pmatrix}+
\lambda\begin{pmatrix}
0 &0 &0\\
0 &0 &1\\
0 &0 &0
\end{pmatrix}
\right\vert\equiv0$, for any $\lambda\in \mathbb{C}$. Consequently, the strongly regular is no longer the necessary condition for the existence of a unique solution to the singular SDE, which is distinguished from the deterministic counterpart and the case of $C=\mathbf{0}$.
\end{example}

In the following, we investigate the sufficient and necessary conditions for the well-posedness of stochastic strongly regular systems.
\begin{theorem}\label{theorem1}
If the system \eqref{SDE} is strongly regular, it has a unique solution $x(\cdot)\in L^2_{\mathcal{F}_t}(0,T;\mathbb{R}^n)$ for any $u(\cdot)\in L^2_{\mathcal{F}_t}(0,T;\mathbb{R}^r)$
if and only if there exist invertible matrices $M$ and $N\in \mathbb{R}^{n\times n}$ such that
\begin{equation}\label{t20}
\begin{aligned}
&MEN=diag(I_{h},G), \quad MAN=diag(A_1,I_{n-h}),\\
&MB=\begin{pmatrix}
B_1\\
B_2
\end{pmatrix},\quad
MCN=\begin{pmatrix}
C_{11}        &C_{12}\\
\mathbf{0}_{(n-h)\times h}    &C_{22}
\end{pmatrix},\quad
MD=\begin{pmatrix}
D_1\\
D_2
\end{pmatrix},
\end{aligned}
\end{equation}
where $C_{22}B_2=D_2$,
$B_2=\big(B_{21}^\mathrm{T},
B_{22}^\mathrm{T},
\cdots,
B_{2l}^\mathrm{T}\big)^\mathrm{T}$
and
\begin{equation}\label{t25}
B_{2k}=
\begin{pmatrix}
b_{k1}\\
0\\
\vdots\\
0
\end{pmatrix}\in\mathbb{R}^{m_k\times r},\ 1\leq k \leq l,
\end{equation}
$G \in \mathbb{R}^{(n-h)\times (n-h)}$ is a nilpotent matrix of form \eqref{matrix3}, $A_1, C_{11}\in\mathbb{R}^{h\times h},\ C_{12} \in \mathbb{R}^{h\times (n-h)}, \ B_1,D_1\in\mathbb{R}^{h\times r},\ C_{22} \in \mathbb{R}
^{(n-h)\times (n-h)},\
B_2, D_2\in\mathbb{R}^{(n-h)\times r}$ and $b_{k1}\in \mathbb{R}^{1\times r},\
1\leq k\leq l$. In addition, $\sum_{k=1}^lm_k=n-h$.
Furthermore, the initial value satisfies
\begin{equation}\label{t27}
\begin{aligned}
\begin{pmatrix}
\mathbf{0}_{(n-h)\times h}
&I_{n-h}
\end{pmatrix}Mx_0
=\mathbf{0}_{(n-h)\times 1}.
\end{aligned}
\end{equation}

\end{theorem}

\noindent{\bf Proof}\quad
We assume that
$
MCN=\begin{pmatrix}
C_{11} &C_{12}\\
C_{21} &C_{22}
\end{pmatrix}$. Let $
x=N\begin{pmatrix}
x_1\\
x_2
\end{pmatrix}
$.
Substituting it into equation \eqref{SDE} and multiplying by the invertible matrix $M$ on both sides of equation \eqref{SDE}, we have
\begin{equation}\label{t3}
\left\{
\begin{aligned}
&dMEN\begin{pmatrix}
x_1\\
x_2
\end{pmatrix}
=\left[MAN\begin{pmatrix}
x_1\\
x_2
\end{pmatrix}
+MBu\right]dt
+\left[MCN\begin{pmatrix}
x_1\\
x_2
\end{pmatrix}+MDu\right]dW(t),\\
&MEN\begin{pmatrix}
x_1(0)\\
x_2(0)
\end{pmatrix}=Mx_0.
\end{aligned}
\right.
\end{equation}
Since the system \eqref{SDE} is strongly regular, equation \eqref{t3} can be rewritten as
\begin{subnumcases}{\label{SDE13}}
dx_1(t)=[A_1x_1(t)+B_1u(t)]dt
+[C_{11}x_1(t)+C_{12}x_2(t)+D_1u(t)]dW(t),
\label{EQ13_1}\\
dGx_2(t)=[x_2(t)+B_2u(t)]dt
+[C_{21}x_1(t)+C_{22}x_2(t)+D_2u(t)]dW(t)
,\label{EQ13_2}\\
x_1(0)=\begin{pmatrix}
I_h &\mathbf{0}_{h\times(n-h)}
\end{pmatrix}Mx_0:=x_{10},\
Gx_2(0)=\begin{pmatrix}
\mathbf{0}_{(n-h)\times h} & I_{n-h}\end{pmatrix}Mx_0:=x_{20}.\notag
\end{subnumcases}
Accounting for $G=diag\big(G_1,G_2,\cdots,G_l\big)$, we decompose $x_2,B_2,C_{21},C_{22}$ and $D_2$ as follows
\begin{equation*}
\begin{aligned}
x_2=\begin{pmatrix}
x_{21}\\
x_{22}\\
\vdots\\
x_{2l}
\end{pmatrix},
B_2=\begin{pmatrix}
B_{21}\\
B_{22}\\
\vdots\\
B_{2l}
\end{pmatrix},
C_{21}=\begin{pmatrix}
C_{21}^1\\
C_{21}^2\\
\vdots\\
C_{21}^l
\end{pmatrix},
C_{22}=\begin{pmatrix}
C_{22}^1\\
C_{22}^2\\
\vdots\\
C_{22}^l
\end{pmatrix},
D_2=\begin{pmatrix}
D_{21}\\
D_{22}\\
\vdots\\
D_{2l}
\end{pmatrix}.
\end{aligned}
\end{equation*}
Then, equation \eqref{EQ13_2} can be written as
\begin{equation*}
dG_{k}x_{2k}=\big[x_{2k}(t)+B_{2k}u(t)
\big]dt+\big[C_{21}^kx_1(t)
+C_{22}^kx_2(t)+D_{2k}u(t)\big]dW(t),\
1\leq k \leq l.
\end{equation*}
We decompose $x_{2k},\ C^k_{21},\ C^k_{22}$ and $D_{2k}$ in the following manner
$$x_{2k}=\begin{pmatrix}
x_{2k,1}\\
x_{2k,2}\\
\vdots\\
x_{2k,m_k-1}\\
x_{2k,m_k}
\end{pmatrix},\ C_{21}^k=\begin{pmatrix}
C_{21}^{k,1}\\
C_{21}^{k,2}\\
\vdots\\
C_{21}^{k,m_k-1}\\
C_{21}^{k,m_k}
\end{pmatrix},\
C_{22}^k=\begin{pmatrix}
C_{22}^{k,1}\\
C_{22}^{k,2}\\
\vdots\\
C_{22}^{k,m_k-1}\\
C_{22}^{k,m_k}
\end{pmatrix},\
D_{2k}=\begin{pmatrix}
D_{2k,1}\\
D_{2k,2}\\
\vdots\\
D_{2k,m_k-1}\\
D_{2k,m_k}
\end{pmatrix}.$$
By the strongly regular property of the system \eqref{SDE}, \eqref{matrix3} and \eqref{t25},
we have for $1\leq k\leq l$,
\begin{subnumcases}{\label{SDE21}}
dx_{2k,2}(t)=[x_{2k,1}(t)+b_{k1}u(t)]dt
+[C_{21}^{k,1}x_1(t)
+C_{22}^{k,1}x_2(t)+D_{2k,1}u(t)]dW(t),\label{EQ21_1}\\
dx_{2k,3}(t)=x_{2k,2}(t)dt
+[C_{21}^{k,2}x_1(t)
+C_{22}^{k,2}x_2(t)+D_{2k,2}u(t)]dW(t),
\label{EQ21_2}\\
\qquad \qquad \vdots \notag\\
dx_{2k,m_k}(t)=x_{2k,m_k-1}(t)dt
+[C_{21}^{k,m_k-1}x_1(t)
+C_{22}^{k,m_k-1}x_2(t)
+D_{2k,m_k-1}u(t)]dW(t), \label{EQ21_3}\\
0=x_{2k,m_k}(t)dt
+[C_{21}^{k,m_k}x_1(t)
+C_{22}^{k,m_k}x_2(t)
+D_{2k,m_k}u(t)]dW(t), \label{EQ21_4}
\end{subnumcases}
Based on the property of continuous martingale, we solve the above equation from \eqref{EQ21_4} to
\eqref{EQ21_1} and we have for $k=1,2,\cdots,l$,
\begin{equation}\label{t28}
\left\{
\begin{aligned}
x_{2k,1}&=-b_{k1}u,\ \text{a.e.}\ t,\\
x_{2k,2}&=\cdots=x_{2k,m_k}\equiv0,
\end{aligned}
\right.
\end{equation}
it follows that $Gx_2=\mathbf{0}_{(n-h)\times1}$. By \eqref{EQ13_2}, we get $x_2(t)=-B_2u(t)$, $C_{21}x_1(t)+C_{22}x_2(t)
+D_2u(t)=\mathbf{0}_{(n-h)\times1}$, for a.e. $t$.
Then, equation \eqref{SDE13} is transformed to
\begin{equation}\label{SDE16}
\left\{\begin{aligned}
&dx_1(t)=[A_1x_1(t)+B_1u(t)]dt
+[C_{11}x_1(t)+\tilde{D}_1u(t)]dW(t),
\quad x_1(0)=x_{10},\\
&x_2(t)=-B_2u(t),\ \text{a.e.}\ t,\ Gx_2(0)=x_{20}=\mathbf{0}_{(n-h)\times1},\\
&C_{21}x_1(t)+C_{22}x_2(t)+D_2u(t)
=\mathbf{0}_{(n-h)\times1},\ \text{a.e.}\ t.
\end{aligned}
\right.
\end{equation}
where $\tilde{D}_1=D_1-C_{12}B_2$.
By solving the first equation of \eqref{SDE16} which is a linear SDE, we obtain
\begin{subnumcases}{\label{SDE22}}
x_1(t)=L(t)^{-1}\Big\{x_{10}+
\int_0^t\big(B_1-C_{11}
\tilde{D}_1\big)u(s)ds +\int_0^t\tilde{D}_1
u(s)L(s)dW(s)\Big\},\label{EQ22_1}\\
x_2(t)=-B_2u(t),\ \text{a.e.}\ t,\  x_{20}=\mathbf{0}_{(n-h)\times1},\label{EQ22_2}\\
C_{21}x_1(t)=(C_{22}B_2-D_2)u(t),\ \text{a.e.}\ t
\label{EQ22_3},
\end{subnumcases}
where
$L(t)=\exp\Big\{\big(\frac{1}{2}
C_{11}^\mathrm{T}C_{11}-A_1\big)t-
C_{11}W(t)\Big\}.$
Substituting \eqref{EQ22_1} into \eqref{EQ22_3}, since it holds for any $u(\cdot)\in L^2_{\mathcal{F}_t}
(0,T;\mathbb{R}^r)$, we have
$$C_{21}=\mathbf{0}_{(n-h)\times h}, \ C_{22}B_2=D_2.$$
Therefore, we have led to the necessity of Theorem \ref{theorem1}, and the sufficiency is easy to verify.
\hfill$\Box$
\begin{remark}
From \eqref{t28}, the solution to equation \eqref{SDE} is not necessary continuous because $x_{2k,1}(t)=-b_{k1}u(t)$,
a.e. $t,\ k=1,2,\cdots,l$. Indeed, it can be seen from Example \ref{ex1} where $x_1$ and $x_2$ are not continuous, respectively. Therefore, the solution to singular SDE does not belong to $L^2_{\mathcal{F}_t}(\Omega;
C([0,T];\mathbb{R}^n))$, which is distinguished from the normal SDE.
\end{remark}
\begin{remark}
We now compare our conclusion with Lemma 2.1 in \cite{zhang2017linear}. Consider the following singular SDE
\begin{equation*}
\left\{
\begin{aligned}
&d\begin{pmatrix}
1 &0 &0\\
0 &1 &0\\
0 &0 &0
\end{pmatrix}
\begin{pmatrix}
x_1\\
x_{21}\\
x_{22}
\end{pmatrix}=\left[
I_{3}\begin{pmatrix}
x_1\\
x_{21}\\
x_{22}
\end{pmatrix}+
\begin{pmatrix}
1 \\
1 \\
0
\end{pmatrix}u\right]dt
+\left[
\begin{pmatrix}
1 &0 &0\\
0 &0 &0\\
0 &0 &1
\end{pmatrix}
\begin{pmatrix}
x_1\\
x_{21}\\
x_{22}
\end{pmatrix}
\right] dW(t),\\
&x_1(0)=0,\
x_{21}(0)=0,\
x_{22}(0)=0.
\end{aligned}
\right.
\end{equation*}
Noted that rank $(E)=2$ and rank $(E,C)=3$, which yields that  Assumption 2.1 in \cite{zhang2017linear} does not hold.
However, let $M=N=I_3$, one can check that $G=0$, $A_1=I_2$, $B_1=\begin{pmatrix}
1\\
1
\end{pmatrix}$, $B_2=0$, $C_{11}=
\begin{pmatrix}
1 &0\\
0 &0
\end{pmatrix}$, $C_{12}=\mathbf{0}_{2\times1}$, $C_{22}=1$, $D=\mathbf{0}_{3\times1}$, which means that the conditions in Theorem 3.2 are satisfied. Thus, the above singular SDE should have a unique solution for any $u\in L^2_{\mathcal{F}_t}(0,T;\mathbb{R})$. To check this, one may note that the above equation is equivalent to
\begin{equation*}
\left\{\begin{aligned}
&dx_1=(x_1+u)dt+x_1dW(t),\\
&dx_{21}=(x_{21}+u)dt,\\
&0=x_{22}dt+x_{22}dW(t),\\
&x_1(0)=0,\ x_{21}(0)=0,\ x_{22}(0)=0.
\end{aligned}
\right.
\end{equation*}
One can easily check it has a unique solution for any $u\in L^2_{\mathcal{F}_t}(0,T;\mathbb{R})$. This indicates that our conditions are more general. If in addition to the conditions in Theorem 3.2, one also assume that  $C_{22}=\mathbf{0}_{(n-h)\times(n-h)},\
D=\mathbf{0}_{n\times r}$,  $m_k=1$ for any $1\leq k\leq l$, $l=n-h$,  $B_2=(b_{11},b_{21},\cdots,b_{l1})^\mathrm{T}$ is any vector in $\mathbb{R}^{(n-h)\times r}$, then it will correspond to Assumption 2.1 and Assumption 2.2 in \cite{zhang2017linear}.
\end{remark}
\begin{remark}\label{remark1}
From the above derivation and Example \ref{ex1}, it is noted that conditions \eqref{t20}-\eqref{t27} are only sufficient conditions for the existence and uniqueness of the solution to system \eqref{SDE}. However, these conditions can be necessary in some special cases, for example, when $C=A$.
\end{remark}
Let us now consider the case of $C=A$.
According to \eqref{kronecker1}-\eqref{matrix3}, the solvability of equation \eqref{SDE} is equivalent to that of the following equations
\begin{subequations}\label{AC1}
 \begin{align}
  &\mathbf{0}_{\tilde{n}_0\times n_0}d\tilde{x}_{1}(t)
  =[\mathbf{0}_{\tilde{n}_0\times n_0}\tilde{x}_{1}(t)
  +B_{1}u(t)]dt
  +[\mathbf{0}_{\tilde{n}_0\times n_0}\tilde{x}_{1}(t)+D_{1}u(t)]dW(t), \label{SDE23a}\\
  &L_id\tilde{x}_{2i}(t)
  =[J_i\tilde{x}_{2i}(t)+B_{2i}u(t)]dt
  +[J_i\tilde{x}_{2i}(t)+D_{2i}u(t)]dW(t),
  \quad i=1,2,\cdots,p ,\label{SDE23b}\\
  &\check{L}_jd\tilde{x}_{3j}(t)
  =[\check{J}_j\tilde{x}_{3j}(t)
  +B_{3j}u(t)]dt+
   [\check{J}_j\tilde{x}_{3j}(t)
   +D_{3j}u(t)]dW(t),
  \quad j=1,2,\cdots,q, \label{SDE23c}\\
  &d\tilde{x}_4(t)
  =[A_1\tilde{x}_4(t)+B_4u(t)]dt+
   [A_1\tilde{x}_4(t)+D_4u(t)]dW(t), \label{SDE23d}\\
  &G_kd\tilde{x}_{5k}(t)
  =[\tilde{x}_{5k}(t)
  +B_{5k}u(t)]dt+
  [\tilde{x}_{5k}(t)+D_{5k}u(t)]dW(t), \quad k=1,2,\cdots,l. \label{SDE23e}
  \end{align}
\end{subequations}
The solvability of equations \eqref{SDE23a} and \eqref{SDE23b} is similar to that of equations \eqref{SDE5a} and \eqref{SDE5b}, respectively. And equation \eqref{SDE23c} can be rewritten as
\begin{subnumcases}{\label{SDE24}}
d\tilde{x}_{3j,1}(t)=B_{3j,1}u(t)dt+D_{3j,1}u(t)dW(t), \label{SDE24-1}\\
d\tilde{x}_{3j,2}(t)=[\tilde{x}_{3j,1}(t)+B_{3j,2}u(t)]dt
+[\tilde{x}_{3j,1} (t)+D_{3j,2}u(t)]dW(t),
\label{SDE24-2}\\
\qquad \qquad \vdots \notag\\
d\tilde{x}_{3j,n_j}(t)=[\tilde{x}_{3j,n_j-1}
(t)+B_{3j, n_j}u(t)]dt
+[\tilde{x}_{3j,n_j-1}(t)+D_{3j,n_j}u(t)]dW(t),
\label{SDE24-3}\\
0=[\tilde{x}_{3j,n_j}(t)+B_{3j, n_j+1}u(t)]dt
+[\tilde{x}_{3j,n_j}(t)+D_{3j,n_j+1}u(t)]dW(t).
\label{SDE24-4}
\end{subnumcases}
From \eqref{SDE24-4} and the arbitrariness of the control, we have
\begin{equation*}
\left\{
\begin{aligned}
&B_{3j,n_j+1}=D_{3j,n_j+1},\\
&\tilde{x}_{3j,n_j}(t)=-B_{3j,n_j+1}u(t),\ \text{a.e.} \ t.
\end{aligned}
\right.
\end{equation*}
By the matrix symmetric technique, we define $\Phi_{3j,n_j+1}=B_{3j,n_j+1}^\mathrm{T}
B_{3j,n_j+1}$.
\begin{itemize}
  \item[1)] $B_{3j,n_j+1}^\mathrm{T}
      \neq\mathbf{0}_{r\times1}$, and then there exist invertible matrices $P_{3j,n_j+1}$ and $Q_{3j,n_j+1}\in\mathbb{R}^{r\times r}$ such that
 \begin{equation*}
      \begin{aligned}
      d\begin{pmatrix}
      1   &0 &\cdots &0  \\
      0  &0 &\cdots &0 \\
      \vdots &\vdots   & &\vdots\\
      0 &0 &\cdots &0
      \end{pmatrix}\tilde{u}(t)
      =&-P_{3j,n_{j}+1}^{-1}B_{3j,
      n_j+1}^\mathrm{T}\big[\tilde{x}_{3j,n_j-1}(t)
      +B_{3j,n_j}Q_{3j,n_{j}+1}^{-1}\tilde{u}(t)\big]dt\\
      &-P_{3j,n_{j}+1}^{-1}B^\mathrm{T}_{
      3j,n_j+1}[\tilde{x}_{3j,n_j-1}(t)+D_{3j,
      n_j}Q_{3j,n_{j}+1}^{-1}\tilde{u}(t)]dW(t),
      \end{aligned}
      \end{equation*}
where $\tilde{u}=Q_{3j,n_j+1}u$.
When $\tilde{u}(t)=
(W(t),0,\cdots,0)^\mathrm{T}$, we have
\begin{equation}\label{SDE25}
\left\{
\begin{aligned}
&\tilde{x}_{3j,n_j-1}(t)=
-B_{3j,n_j}Q^{-1}_{3j,n_j+1}\tilde{u}
(t),\ \text{a.e.}\ t, \\
&P^{-1}_{3j,n_j+1}B_{3j,n_j+1}^\mathrm{T}
(B_{3j,n_j}-D_{3j,n_j})Q_{3j,n_j+1}^{-1}
(W(t),0,\cdots,0)^\mathrm{T}=
(1,0,\cdots,0)^\mathrm{T}.
\end{aligned}
\right.
\end{equation}
It is noted that the left-hand side of \eqref{SDE25} is a stochastic process, while the right-hand side is deterministic. Thus, equation \eqref{SDE25} does not hold. Therefore, equation \eqref{SDE23c} does not have a solution for any $u(\cdot)\in L^2_{\mathcal{F}_t}
(0,T;\mathbb{R}^r)$ when $B_{3j,n_j+1}^\mathrm{T}
\neq\mathbf{0}_{r\times 1}$.
  \item[2)] $B_{3j,n_j+1}^\mathrm{T}=\mathbf{0}_{r\times1}$, and then, we obtain
      \begin{equation*}
\left\{
\begin{aligned}
&B_{3j,n_j}=D_{3j,n_j},\\
&\tilde{x}_{3j,n_j-1}(t)=-B_{3j,n_j}u(t),\ \text{a.e.} \ t.
\end{aligned}
\right.
\end{equation*}
\end{itemize}
Iterating the above procedures, equation \eqref{SDE23c} has a unique solution $\tilde{x}_{3j}$ for any $u(\cdot)\in L^2_{\mathcal{F}_t}(0,T;\mathbb{R}^r)$ only when $B_{3j}=D_{3j}=\mathbf{0}_{(n_j+1)\times r}$, $j=1,2,\cdots,q$. Moreover, the unique solution is $\tilde{x}_{3j}=\mathbf{0}_{n_j\times1}
$.

Analogue to the solvability of equation \eqref{SDE5e}, equation \eqref{SDE23e} has a unique solution for any $u(\cdot)\in L^2_{\mathcal{F}_t}
(0,T;\mathbb{R}^r)$ only when $B_{5k}=D_{5k}=(b_{k1},0,\cdots,0)^\mathrm{T}$, for $k=1,2,\cdots,l$.
Moreover, the unique solution is $(-b_{k1}u,0,\cdots,0)^\mathrm{T}$. 
Therefore,
we have the following necessary and sufficient conditions for the well-posedness of equation \eqref{SDE} when $C=A$.
\begin{theorem}\label{theorem7}
When $C=A$, the equation \eqref{SDE} has a unique solution $x(\cdot)\in L^2_{\mathcal{F}t}(0,T;\mathbb{R}^n)$ for any $u(\cdot)\in L^2_{\mathcal{F}_t}(0,T;\mathbb{R}^r)$ if and only if there exist invertible matrices $M$, $N\in\mathbb{R}^{n\times n}$ such that
\begin{equation*}
\begin{aligned}
&MEN=diag(I_{h},G), \quad MAN=diag(A_1,I_{n-h}),\\
&MB=\begin{pmatrix}
B_1\\
B_2
\end{pmatrix},\quad
MD=\begin{pmatrix}
D_1\\
B_2
\end{pmatrix},
\end{aligned}
\end{equation*}
where
$B_2=\big(B_{21}^\mathrm{T},
B_{22}^\mathrm{T},
\cdots,
B_{2l}^\mathrm{T}\big)^\mathrm{T}$,
$B_{2k}=
\big(b_{k1},0,\cdots,0
\big)^\mathrm{T},\ 1\leq k \leq l,$
$G \in \mathbb{R}^{(n-h)\times (n-h)}$ is a nilpotent matrix of form \eqref{matrix3},
$A_1\in\mathbb{R}^{h\times h}$ is any matrix,
$B_1,D_1\in\mathbb{R}^{h\times r},
B_2\in\mathbb{R}^{(n-h)\times r}$, $B_{2k}\in\mathbb{R}^{m_k\times r}$ and $b_{k1}\in \mathbb{R}^{1\times r},\
1\leq k\leq l$. In addition,
$\sum_{k=1}^lm_k=n-h$.
Moreover, the initial value satisfies
$\begin{pmatrix}
\mathbf{0}_{(n-h)\times h}
&I_{n-h}
\end{pmatrix}Mx_0
=\mathbf{0}_{(n-h)\times 1}.$
\end{theorem}
\begin{remark}
Equation \eqref{SDE13} (equivalent to \eqref{SDE16}) clearly reveals the structure of singular systems. The first one is a SDE which constitutes the dynamic part of the system. The second one (see \eqref{SDE16}) is the static part which is represented by an algebraic equation and reflects the relationship between sub-variables. Therefore, singular systems can be regarded as a class of composite systems made up by some coupled subsystems.
\end{remark}
\begin{remark}
It looks a challenge to solve completely the general case of $C\neq\mathbf{0}$. We have made various attempts. The following are some heuristic arguments. It is noted that the explicit solution of normal SDE can be obtained by dual homogeneous SDE. Therefore, we introduce the following SDE
\begin{equation*}
\left\{\begin{aligned}
dy(t)=&\mathcal{A}y(t)dt +\mathcal{C}y(t)dW(t),\\
y(0)=&I_n,\\
\end{aligned}
\right.
\end{equation*}
where $\mathcal{A}$ and $\mathcal{C}$
are undetermined coefficients.
Applying It\^o formula to $y^\mathrm{T}(\cdot)Ex(\cdot)$, we have
\begin{equation}\label{T3}
\begin{aligned}
dy^\mathrm{T}(t)Ex(t)=&\Big[
\big(y^\mathrm{T}(t)\mathcal{A}^\mathrm{T}
E+y^\mathrm{T}(t)A+
y^\mathrm{T}(t)
\mathcal{C}^\mathrm{T}C\big)x(t)
+\big(y^\mathrm{T}(t)B
+y^\mathrm{T}\mathcal{C}^\mathrm{T}D\big)u(t)\Big]dt\\
&+\Big[\big(y^\mathrm{T}(t)
\mathcal{C}^\mathrm{T}E
+y^\mathrm{T}(t)C\big)x(t)
+y^\mathrm{T}(t)Du(t)\Big]dW(t).
\end{aligned}
\end{equation}
It is noted that $y(t)=\exp\big(\int_0^t(\mathcal{A}-
\frac{1}{2}\mathcal{C}^\mathrm{T}
\mathcal{C})ds+\int_0^t\mathcal{C}dW(s)
\big)$. Therefore, when $rank\ (E)=rank\begin{pmatrix}
E\\
C
\end{pmatrix}$, we can choose $\mathcal{C}$ satisfying
$y^\mathrm{T}(t)\mathcal{C}^\mathrm{T}E
=-y^\mathrm{T}(t)C$ such that the diffusion term of \eqref{T3} is independent of $x(\cdot)$.
And then, $x(\cdot)$ satisfies a new singular SDE
\begin{equation*}
\left\{\begin{aligned}
d\hat{E}(t)x(t)=&[\hat{A}(t)x(t)
+\hat{B}(t)u(t)]dt +\hat{D}(t)u(t)dW(t),\\
\hat{E}(0)x(0)=&x_0,\\
\end{aligned}
\right.
\end{equation*}
where
$
\hat{E}(t)=y^\mathrm{T}(t)E,\
\hat{A}(t)=y^\mathrm{T}(t)\mathcal{A}^\mathrm{T}
E+y^\mathrm{T}(t)A+
y^\mathrm{T}(t)
\mathcal{C}^\mathrm{T}C,\
\hat{B}(t)=y^\mathrm{T}(t)B
+y^\mathrm{T}\mathcal{C}^\mathrm{T}D,\
\hat{D}(t)=y^\mathrm{T}(t)D.
$
Unfortunately, the coefficients $\hat{E}$, $\hat{A}$, $\hat{B}$ and $\hat{D}$ depend on $t$ and $\omega$, which yields that the previous derivation of proving Theorem \ref{C0} cannot be applied to obtain the conclusion.
\end{remark}
\begin{remark}
We are dedicated to establishing the conditions that ensure the uniqueness of the singular equations. We mention that it is crucial to ensuring the uniqueness of the equations, since it will  guarantee the uniqueness of the optimal solution to the control problem. Otherwise, there would be no one-to-one correspondence between control and state variables, and, as a result, optimal feedback control could not be obtained.
\end{remark}
\section{Finite-horizon stochastic singular LQ control problem}
We consider the following stochastic singular system
\begin{equation}\label{SDE6}
\left\{\begin{aligned}
&Edx(t) = [Ax(t)+Bu(t)]dt
+[Cx(t)+Du(t)]dW(t),\ t\in(0,T],\\
&Ex(0) = x_0,
\end{aligned}
\right.
\end{equation}
where $E\in \mathbb{R}^{n\times n}$ with rank $(E)=m<n$, $A$ and $C \in\mathbb{R}^{n\times n}$, $B$ and $D \in\mathbb{R}^{n\times r}$ are constant matrices, and $x_0\in\mathbb{R}^n$ is initial condition.

Define the set of all admissible control-state pairs
\begin{equation}\label{addmiss1}
\begin{aligned}
\mathcal{J}_{ad}[0,T]=& \Big\{\big(u(\cdot),x(\cdot)\big)
\Big{\vert} \ u(\cdot)\in L_{\mathcal{F}_t}^2(0,T;\mathbb{R}^r)
 \ \text{such that equation}\  \eqref{SDE6}\
  \text{admits a}\\
   &\text{unique solution}\ x(\cdot)\in L_{\mathcal{F}_t}^2(0,T;\mathbb{R}^n)
 \ \text{under certain initial conditions}\Big\}.
\end{aligned}
\end{equation}

For $Q,\ H \in\mathbb{R}^{n\times n}$ and $R\in\mathbb{R}^{r\times r}$, the cost functional is given by
\begin{equation}\label{costf}
\begin{aligned}
J\big(u(\cdot),x(\cdot)\big)=\frac{1}{2}
\mathbb{E}\big[\int_0^T\big(x^\mathrm{T}(t)Qx(t)
+u^\mathrm{T}(t)Ru(t)\big) dt\big]+\frac{1}{2}\mathbb{E}\big[
x^\mathrm{T}(T)E^\mathrm{T}HEx(T)\big].
\end{aligned}
\end{equation}
The inclusion of $E$ in the terminal term of the cost functional \eqref{costf} is indeed a technical requirement to facilitate the transformation into a normal LQ control problem. This is necessary for finite-horizon singular LQ control problems, see also \cite{bender1987linear,feng2003singular, zhang2017linear}.
Let us now state our finite-horizon stochastic singular LQ control problem.

\textbf{Problem (FSP)}. For given $x_0\in\mathbb{R}^n$, find the optimal pair $\big(\bar{u}(\cdot),\bar{x}(\cdot)\big)
\in\mathcal{J}_{ad}[0;T]$ such that
$$J\big(\bar{u}(\cdot),\bar{x}(\cdot)\big)=\inf_
{(u(\cdot),x(\cdot))
\in\mathcal{J}_{ad}[0,T]}
J\big(u(\cdot),x(\cdot)\big)
:=V(0,x_0).
$$

It should be pointed out that the set of admissible control-state pairs may be empty according to the analysis in Section 3. To ensure that  $\mathcal{J}_{ad}[0,T]$ is nonempty, we need the system \eqref{SDE6} to meet some conditions.
Stimulated by \cite{dai1989singular} and \cite{duan2010analysis}, in order to weaken assumptions of the system, we give a state feedback control as follows: for any $v(\cdot)\in
L^2_{\mathcal{F}_t}(0,T;\mathbb{R}^r)$ and $K\in\mathbb{R}^{r\times n}$, define
\begin{equation}\label{con1}
u(\cdot)=Kx(\cdot)+v(\cdot),
\end{equation}
where $x(\cdot)$ satisfies the corresponding state equation
\begin{equation}\label{SDE7}
\left\{\begin{aligned}
&Edx(t) = [(A+BK)x(t)+Bv(t)]dt
+[(C+DK)x(t)+Dv(t)]dW(t),\\
&Ex(0) = x_0.
\end{aligned}
\right.
\end{equation}
We choose $K$ such that the system \eqref{SDE7} satisfies \eqref{t20}-\eqref{t27}, which means that there exist $M_1$ and $N_1$ such that
\begin{equation}\label{t32}
\begin{aligned}
&M_1EN_1=diag(I_{h},G), \quad M_1(A+BK)N_1=diag(A_1,I_{n-h}),\\
&M_1B=\begin{pmatrix}
B_1\\
B_2
\end{pmatrix},\quad
M_1(C+DK)N_1=\begin{pmatrix}
C_{11}        &C_{12}\\
\mathbf{0}_{(n-h)\times h}    &C_{22}
\end{pmatrix},\quad
M_1D=\begin{pmatrix}
D_1\\
D_2
\end{pmatrix},
\end{aligned}
\end{equation}
where $C_{22}B_2=D_2$,
$B_2=\big(B_{21}^\mathrm{T},
B_{22}^\mathrm{T},
\cdots,
B_{2l}^\mathrm{T}\big)^\mathrm{T}$
and
\begin{equation}\label{t33}
B_{2k}=
\begin{pmatrix}
b_{k1}\\
0\\
\vdots\\
0
\end{pmatrix}\in\mathbb{R}^{m_k\times r},\ 1\leq k \leq l,
\end{equation}
$G \in \mathbb{R}^{(n-h)\times (n-h)}$ is a nilpotent matrix of form \eqref{matrix3}, $A_1, C_{11}\in\mathbb{R}^{h\times h},\ C_{12} \in \mathbb{R}^{h\times (n-h)}, \ B_1,D_1\in\mathbb{R}^{h\times r},\ C_{22} \in \mathbb{R}
^{(n-h)\times (n-h)},\
B_2, D_2\in\mathbb{R}^{(n-h)\times r}$ and $b_{k1}\in \mathbb{R}^{1\times r},\
1\leq k\leq l$. In addition, $\sum_{k=1}^lm_k=n-h$.
Furthermore, the initial value satisfies:
\begin{equation}\label{t34}
\begin{aligned}
\begin{pmatrix}
\mathbf{0}_{(n-h)\times h}
&I_{n-h}
\end{pmatrix}M_1x_0
=\mathbf{0}_{(n-h)\times 1}.
\end{aligned}
\end{equation}
It is noted that the system \eqref{SDE6} itself may not satisfy conditions \eqref{t20}-\eqref{t27}, but it may happen that there exists a matrix $K$ such that \eqref{t32}-\eqref{t34} hold, see Example \ref{example1}. Once $K$ is chosen, according to Theorem \ref{theorem1}, the equation \eqref{SDE7} admits a unique solution for any $v(\cdot)\in L^2_{\mathcal{F}_t}(0,T;\mathbb{R}^r).$ In particular, when $K=\mathbf{0}_{r\times n}$ and the system \eqref{SDE6} is strongly regular, it satisfies the conditions \eqref{t32}-\eqref{t34}.

For stochastic normal systems, the optimal control can be judged from the well-known stochastic Pontryagin's maximum principle. However, for stochastic singular systems, we can not get the estimations of state, so the classical maximum principle cannot be used directly to deal with this problem. Motivated by the method in deterministic case \cite{
bender1987linear,
zhaolin1988optimal,
cobb1983descriptor,
dai1989singular,
duan2010analysis,
zhu2002singular}, we will change the stochastic singular LQ control problem into normal LQ control problem, and then we can solve it by maximum principle.

By denoting
\begin{equation}\label{t7}
x=N_1\begin{pmatrix}
x_1\\
x_2
\end{pmatrix},\quad
x_1\in\mathbb{R}^h,\quad
x_2\in\mathbb{R}^{n-h},
\end{equation}
then similar to \eqref{SDE13} and \eqref{SDE16}, the system \eqref{SDE7} is restricted equivalent to the following system:
\begin{subnumcases}{\label{SDE8}}
dx_1(t) = [A_1x_1(t)+B_1v(t)]dt
+[C_{11}x_1(t)+C_{12}x_2(t)+D_1v(t)]dW(t),
\label{EQ8_1}\\
0=x_2(t)+B_2v(t), \label{EQ8_2}\\
x_1(0)=
\begin{pmatrix}
I_h &\mathbf{0}_{h \times(n-h)}
\end{pmatrix}
M_1x_0:=x_{10},\
Gx_2(0)=\begin{pmatrix}
\mathbf{0}_{(n-h)\times h} &I_{n-h}
\end{pmatrix}
M_1x_0:=x_{20}.\notag
\end{subnumcases}
From \eqref{EQ8_2}, we obtain $x_2(\cdot)=-B_2v(\cdot)$. Substituting it into \eqref{EQ8_1}, we have:
\begin{equation}\label{SDE9}
\left\{\begin{aligned}
&dx_1(t) = [A_1x_1(t)+B_1v(t)]dt
+[C_{11}x_1(t)+\tilde{D}_1v(t)]dW(t),\\
&x_1(0)=x_{10},
\end{aligned}
\right.
\end{equation}
where
\begin{equation}
\tilde{D}_1=D_1-C_{12}B_2.
\end{equation}
Combining $x_2=-B_2v$ with \eqref{con1} and \eqref{t7}, we have:
\begin{equation}\label{t9}
\begin{aligned}
\begin{pmatrix}
x\\
u
\end{pmatrix}&=\begin{pmatrix}
I_n         &\mathbf{0}_{n\times r}\\
K           &I_r
\end{pmatrix}
\begin{pmatrix}
x\\
v
\end{pmatrix}
=
\begin{pmatrix}
N_1         &\mathbf{0}_{n\times r}\\
KN_1        &I_r
\end{pmatrix}
\begin{pmatrix}
x_1\\
x_2\\
v
\end{pmatrix}\\
&=
\begin{pmatrix}
N_1         &\mathbf{0}_{n\times r}\\
KN_1        &I_r
\end{pmatrix}
\begin{pmatrix}
I_h           &\mathbf{0}_{h\times r}\\
\mathbf{0}_{(n-h) \times h}  &-B_2\\
\mathbf{0}_{r\times h}  &I_r
\end{pmatrix}
\begin{pmatrix}
x_1\\
v
\end{pmatrix}.
\end{aligned}
\end{equation}
Next, we introduce a cost functional\begin{equation}\label{costf3}
\begin{aligned}
J_1(v(\cdot))=\frac{1}{2}
\mathbb{E}\int_0^T
\begin{pmatrix}
x_1\\
v
\end{pmatrix}^\mathrm{T}
\begin{pmatrix}
\hat{Q}      &S^\mathrm{T}\\
S &\hat{R}
\end{pmatrix}
\begin{pmatrix}
x_1\\
v
\end{pmatrix} dt
+\frac{1}{2}\mathbb{E}\big[
x_1^\mathrm{T}(T)\hat{H}
x_1(T)\big]
,
\end{aligned}
\end{equation}
where
\begin{equation}\label{t8}
\begin{aligned}
\begin{pmatrix}
\hat{Q}      &S^\mathrm{T}\\
S &\hat{R}
\end{pmatrix}=
\begin{pmatrix}
I_h         &\mathbf{0}\\
\mathbf{0}  &-B_2\\
\mathbf{0}  &I_r
\end{pmatrix}^\mathrm{T}
\begin{pmatrix}
N_1         &\mathbf{0}\\
KN_1        &I_r
\end{pmatrix}^\mathrm{T}
\begin{pmatrix}
Q           &\mathbf{0}\\
\mathbf{0}  &R
\end{pmatrix} \begin{pmatrix}
N_1         &\mathbf{0}\\
KN_1        &I_r
\end{pmatrix}
\begin{pmatrix}
I_h         &\mathbf{0}\\
\mathbf{0}  &-B_2\\
\mathbf{0}  &I_r
\end{pmatrix},
\end{aligned}
\end{equation}
\begin{equation*}\label{t12}
\begin{aligned}
\hat{H}=
\begin{pmatrix}
I_h           &\mathbf{0}_{h\times (n-h)}
\end{pmatrix}
(M_1^{-1})^\mathrm{T}HM_1^{-1}\begin{pmatrix}
I_h           \\
\mathbf{0}_{(n-h)\times h}
\end{pmatrix}.
\end{aligned}
\end{equation*}
Define the admissible control set
\begin{equation}\label{admiss2}
\mathcal{U}_{ad}[0,T]=\big\{v(\cdot)\
\vert \ v(\cdot)\in L_{\mathcal{F}_t}^2(0,T;\mathbb{R}^r)
 \big\}.
\end{equation}
We state the stochastic normal LQ control problem as follows.

\textbf{Problem $\rm(FP)_K$.} For given $x_{10}\in\mathbb{R}^h$, find the optimal control $\bar{v}(\cdot)\in\mathcal{U}_{ad}[0,T]$ such that
$$J_1(\bar{v}(\cdot))=\inf_{v(\cdot)\in
\mathcal{U}_{ad}[0,T]}J_1(v(\cdot))
:=V_1(0,x_{10}).$$
\begin{lemma}\label{lem1}
When $K$ is chosen such that \eqref{t32}-\eqref{t34} hold, Problem (FSP) and Problem $(FP)_K$ are equivalent.
\end{lemma}
\noindent{\bf Proof}\quad
From \eqref{t9}, cost functional \eqref{costf} becomes
\begin{equation*}
\begin{aligned}
J\big(u(\cdot),x(\cdot)\big)=&\frac{1}{2}
\mathbb{E}\int_0^T
\begin{pmatrix}
x_1\\
v
\end{pmatrix}^\mathrm{T}
\begin{pmatrix}
\hat{Q}      &S^\mathrm{T}\\
S &\hat{R}
\end{pmatrix}
\begin{pmatrix}
x_1\\
v
\end{pmatrix} dt\\
&+\frac{1}{2}\mathbb{E}
\begin{pmatrix}x_1(T)\\
Gx_2(T)
\end{pmatrix}^\mathrm{T}
(M_1^{-1})^\mathrm{T}HM_1^{-1}
\begin{pmatrix}x_1(T)\\
Gx_2(T)
\end{pmatrix}.
\end{aligned}
\end{equation*}
According to $x_{2k,2}=\cdots=x_{2k,m_k}=0,\ 1\leq k\leq l$ and \eqref{matrix3},
we have $Gx_2(T)=0$. Therefore,
$$\inf_{(u(\cdot),x(\cdot))\in\mathcal{J}_{ad}
[0,T]}
J(u(\cdot),x(\cdot))=\inf_{(u(\cdot),x(\cdot))\in
\mathcal{J}_{ad}[0,T]}J_1
\big(u(\cdot)-Kx(\cdot)\big)
\geq \inf_{v(\cdot)\in\mathcal{U}_{ad}
[0,T]}J_1(v(\cdot)).
$$
On the contrary, we obtain
$$\inf_{v(\cdot)\in\mathcal{U}_{ad}
[0,T]}
J_1(v(\cdot))=\inf_{v(\cdot)\in
\mathcal{U}_{ad}[0,T]}J
\big(KN_1\begin{pmatrix}
x_1\\
x_2
\end{pmatrix}+v(\cdot),N_1\begin{pmatrix}
x_1\\
x_2
\end{pmatrix}\big)
\geq \inf_{(u(\cdot),x(\cdot))\in\mathcal{J}_{ad}
[0,T]}
J(u(\cdot),x(\cdot)),
$$
which yields
\begin{equation*}\inf_{(u(\cdot),x(\cdot))\in\mathcal{J}_{ad}
[0,T]}
J(u(\cdot),x(\cdot))=\inf_{v(\cdot)\in\mathcal{U}_{ad}
[0,T]}
J_1(v(\cdot)).
\end{equation*}
\hfill$\Box$

As a consequence, the stochastic singular LQ control problem is replaced by a stochastic normal LQ control problem with the order $h<n$. Now, we will solve Problem ${\rm (FP)_K}$ by stochastic maximum principle. To begin with, we introduce the following assumptions.
\begin{itemize}
             \item (H1) $
Q>0, \ R>0, \ H\geq0.
$
             \item (H1)$'$ $
Q\geq0, \ \text{rank} \left(I_r+KN_1\begin{pmatrix}
\mathbf{0}_{h\times r}\\
-B_2
\end{pmatrix}\right)=r,\ R>0, \ H\geq0.
$
\end{itemize}
Let $\tilde{Q}
=\hat{Q}-S^\mathrm{T}\hat{R}^{-1}S$, and we have
\begin{equation}\label{t24}
\begin{pmatrix}
\tilde{Q}      &\mathbf{0}\\
\mathbf{0}     &\hat{R}
\end{pmatrix}=
\begin{pmatrix}
I_n                             &\mathbf{0}_{n\times r}\\
-\hat{R}^{-1}S     &I_r
\end{pmatrix}^\mathrm{T}
\begin{pmatrix}
\hat{Q}      &S^\mathrm{T}\\
S &\hat{R}
\end{pmatrix}
\begin{pmatrix}
I_n                             &\mathbf{0}_{n\times r}\\
-\hat{R}^{-1}S     &I_r
\end{pmatrix}.
\end{equation}
Now, we discuss the positive definiteness (semi-definiteness) of $\hat{R}$ and $\tilde{Q}$ in two cases.
On the one hand, the matrix
\begin{equation*}
\begin{pmatrix}
N_1         &\mathbf{0}_{n\times r}\\
KN_1        &I_r
\end{pmatrix}
\end{equation*}
is invertible and the matrix
\begin{equation*}
\begin{pmatrix}
I_h           &\mathbf{0}_{h\times r}\\
\mathbf{0}_{(n-h) \times h}  &-B_2\\
\mathbf{0}_{r\times h}  &I_r
\end{pmatrix}
\end{equation*}
has full column rank. Thus, from (H1), \eqref{t8} and \eqref{t24}, we have $ \hat{R}>0$, $\tilde{Q}>0$. On the other hand, under (H1)$'$, from \eqref{t8} and \eqref{t24}, we have $\tilde{Q}\geq0$ and
$$\hat{R}=\begin{pmatrix}
\mathbf{0} &-B_2^\mathrm{T}\end{pmatrix}N_1^\mathrm{T}
QN_1\begin{pmatrix}
\mathbf{0}\\
-B_2
\end{pmatrix}+\left[
\begin{pmatrix}
\mathbf{0} &-B_2^\mathrm{T}\end{pmatrix}
(KN_1)^\mathrm{T}+I_r\right]R\left[
KN_1\begin{pmatrix}
\mathbf{0}\\
-B_2
\end{pmatrix}+I_r\right],$$
which yields $\hat{R}>0$.

Therefore, from Corollary 5.7 in Chapter 6 of  \cite{yong1999stochastic}, Problem ${\rm (FP)_K}$ is uniquely solvable under (H1) (or (H1)$'$). Moreover, the optimal control $\bar{v}(\cdot)$ can be represented by
\begin{equation*}\label{optimalv}
\bar{v}(t)=-\hat{R}^{-1}[S\bar{x}_1(t)
-B_1^\mathrm{T}p(t)
-\tilde{D}_1^\mathrm{T}q(t)],\quad t\in[0,T],
\end{equation*}
where $(\bar{x}_1(\cdot),p(\cdot),q(\cdot))$ satisfies the following forward backward stochastic differential equation (FBSDE)
\begin{equation}\label{FBSDE}
\left\{\begin{aligned}
d\bar{x}_1(t)=&\big[(A_1-B_1\hat{R}^{-1}S)
\bar{x}_1(t)+B_1\hat{R}^{-1}B_1^\mathrm{T}
p(t)+B_1\hat{R}^{-1}\tilde{D}_1^\mathrm{T}q(t)\big]dt
\\
& +\big[(C_{11}-\tilde{D}_1\hat{R}^{-1}S)
\bar{x}_1(t)+\tilde{D}_1\hat{R}^{-1}B_1^\mathrm{T}
p(t)+\tilde{D}_1\hat{R}^{-1}\tilde{D}_1^\mathrm{T}q(t)\big]dW(t),\\
dp(t)=&
\big[(A_1-B_1\hat{R}^{-1}S)^\mathrm{T}p(t)
+(C_{11}-\tilde{D}_1\hat{R}^{-1}S)^\mathrm{T}
q(t)+(S^\mathrm{T}\hat{R}^{-1}S-\hat{Q})\bar{x}_1(t)
\big]dt\\
&+q(t)dW(t),\quad\\
x_1(0)=&x_{10}, \quad
p(T)=-\hat{H}\bar{x}_1(T).
\end{aligned}
\right.
\end{equation}
Furthermore, FBSDE \eqref{FBSDE} admits a unique solution
$(\bar{x}_1(\cdot),p(\cdot),q(\cdot))\in
L^2_{\mathcal{F}_t}(\Omega;
C([0,T];\mathbb{R}^h))\times
L^2_{\mathcal{F}_t}(\Omega;
C([0,T];\mathbb{R}^h))\times
L_{\mathcal{F}_t}^2(0,T;\mathbb{R}^h)$.
\begin{theorem}
Under {\rm (H1)} (or {\rm (H1)$'$}), when $K$ is picked such that \eqref{t32}-\eqref{t34} hold, Problem (FSP) has a unique optimal control-state pair which can be expressed linearly by the solution to FBSDE \eqref{FBSDE}.
\end{theorem}
\noindent{\bf Proof}\quad
From \eqref{t9}, it follows that the optimal control-state pair $\big(\bar{u}(\cdot),\bar{x}(\cdot)\big)$ of Problem (FSP) is determined by
\begin{equation*}\label{t10}
\begin{aligned}
\begin{pmatrix}
\bar{x}\\
\bar{u}
\end{pmatrix}=
&\begin{pmatrix}
N_1         &\mathbf{0}_{n\times r}\\
KN_1        &I_r
\end{pmatrix}
\begin{pmatrix}
I_h           &\mathbf{0}_{h\times r}\\
\mathbf{0}_{(n-h)\times h}  &-B_2\\
\mathbf{0}_{r\times h}  &I_r
\end{pmatrix}
\begin{pmatrix}
\bar{x}_1\\
\bar{v}
\end{pmatrix}\\
=&\begin{pmatrix}
N_1         &\mathbf{0}_{n\times r}\\
KN_1        &I_r
\end{pmatrix}
\begin{pmatrix}
I_h           &\mathbf{0}_{h\times r}\\
\mathbf{0}_{(n-h)\times h}  &-B_2\\
\mathbf{0}_{r\times h}  &I_r
\end{pmatrix}
\begin{pmatrix}
&I_h           &\mathbf{0}_{h \times h}   &\mathbf{0}_{h \times h}\\
&-\hat{R}^{-1}S &\hat{R}^{-1}B_1^\mathrm{T}
&\hat{R}^{-1}\tilde{D}_1^\mathrm{T}
\end{pmatrix}
\begin{pmatrix}
\bar{x}_1\\
p\\
q
\end{pmatrix},
\end{aligned}
\end{equation*}
where $\big(\bar{x}_1(\cdot), p(\cdot), q(\cdot)\big)$ is given by \eqref{FBSDE}. Therefore, the existence and uniqueness of the optimal control-state pair can be obtained by that of FBSDE \eqref{FBSDE}.
\hfill$\Box$

Subsequently, the optimal control of Problem (FSP) can be written as the state feedback by introducing the following Riccati equation:
\begin{equation}\label{Riccati eq}
\left\{\begin{aligned}
&\dot{P}(t)+P(t)A_1+A_1^\mathrm{T}P(t)
+C_{11}^\mathrm{T}P(t)C_{11}+\hat{Q}\\
&-[B_1^\mathrm{T}P(t)+S+\tilde{D}_1^\mathrm{T}
P(t)C_{11}]^\mathrm{T}[\hat{R}
+\tilde{D}_1^\mathrm{T}P(t)\tilde{D}_1]^{-1}
[B_1^\mathrm{T}P(t)+S+\tilde{D}_1^\mathrm{T}
P(t)C_{11}]=0,\\
&P(T)=\hat{H}.
\end{aligned}
\right.
\end{equation}
\begin{theorem}\label{theorem2}
Let {\rm (H1)} (or {\rm (H1)$'$}) hold, Problem $(FSP)$ admits a unique optimal control-state pair when $K$ is chosen such that \eqref{t32}-\eqref{t34} hold. Furthermore, the optimal control can be represented as a state feedback form
\begin{equation*}\label{optimalv3}
\bar{u}(t)=\big[K-\Psi(t)\begin{pmatrix}
I_h  &\mathbf{0}_{h\times (n-h)}
\end{pmatrix}
N_1^{-1}\big]\bar{x}(t),
\end{equation*}
where
$
\Psi(t)\triangleq[\hat{R}+\tilde{D}_1
^\mathrm{T}P(t)\tilde{D}_1]^{-1}
[B_1^\mathrm{T}P(t)+S
+\tilde{D}_1^\mathrm{T}P(t)C_{11}].
$
Moreover, the minimum value of the cost functional is
\begin{equation*}\label{costf5}
V(0,x_0)=
\frac{1}{2}
x_0^\mathrm{T}M_1^\mathrm{T}
\begin{pmatrix}
I_h\\
\mathbf{0}
\end{pmatrix}P(0)
\begin{pmatrix}
I_h &\mathbf{0}
\end{pmatrix}M_1x_0.
\end{equation*}
\end{theorem}
\noindent{\bf Proof}\quad
From Theorem 6.1 in Chapter 6 of \cite{yong1999stochastic}, the optimal control $\bar{v}(\cdot)$ of Problem ${\rm (FP)_K}$ can be described as a state feedback form as follows
\begin{equation*}\label{optimalv2}
\bar{v}(\cdot)=-\Psi(\cdot)\bar{x}_1(\cdot),
\end{equation*}
and
\begin{equation*}
V_1(0,x_{10})=\frac{1}{2}x_{10}^\mathrm{T}
P(0)x_{10}.
\end{equation*}
Due to \eqref{t7}, we derive $$\bar{u}(t)=K_1\bar{x}(t)+\bar{v}(t)=\big[K_1-\Psi(t)\begin{pmatrix}
I_h  &\mathbf{0}_{h\times (n-h)}
\end{pmatrix}
N_1^{-1}\big]\bar{x}(t),\quad t\in[0,T],$$
and
$$V(0,x_0)=V_1(0,x_{10})=\frac{1}{2}
x_0^\mathrm{T}M_1^\mathrm{T}
\begin{pmatrix}
I_h\\
\mathbf{0}
\end{pmatrix}P(0)
\begin{pmatrix}
I_h &\mathbf{0}
\end{pmatrix}M_1x_0,$$
which completes the proof.
\hfill$\Box$
\section{Infinite-horizon stochastic singular LQ control problem}
In this section, we study the infinite-horizon stochastic singular LQ control problem. The state equation is still given by
\begin{equation}\label{SDE26}
\left\{\begin{aligned}
&Edx(t) = [Ax(t)+Bu(t)]dt
+[Cx(t)+Du(t)]dW(t),\ t>0,\\
&Ex(0) = x_0,
\end{aligned}
\right.
\end{equation}
and the cost functional is formulated as follows
\begin{equation}\label{costf11}
J\big(u(\cdot),x(\cdot)\big)=\mathbb{E}\big[\int_0^\infty
\big(x^\mathrm{T}(t)Qx(t)+u^\mathrm{T}(t)Ru(t)
\big)dt\big].
\end{equation}
Let
\begin{equation*}
\begin{aligned}
\mathcal{J}_{loc}[0,\infty)=\bigcap_{T>0}
\mathcal{J}_{ad}[0,T],
\end{aligned}
\end{equation*}
where $\mathcal{J}_{ad}[0,T]$ is defined by \eqref{addmiss1}.
We introduce the set of admissible control-state pairs
\begin{equation*}
\begin{aligned}
\mathcal{J}_{ad}[0,\infty)=& \Big\{\big(u(\cdot),x(\cdot)\big)\in
\mathcal{J}_{loc}[0,\infty)
\Big{\vert} \ \mathbb{E}\int_0^\infty\vert u(t)\vert ^2dt<\infty,
\text{and the corresponding solution}\\
&\text{to equation}\ \eqref{SDE26}\  \text{satisfies} \ \mathbb{E}\int_0^\infty\vert x(t)\vert ^2dt<\infty\Big\}.
\end{aligned}
\end{equation*}
Our infinite-horizon stochastic singular LQ control problem can be stated as follows.

\textbf{Problem (IFSP)}. For given $x_0\in\mathbb{R}^n$, seek the optimal pair $\big(\bar{u}(\cdot),\bar{x}(\cdot)\big)\in\mathcal{J}_{ad}[0,\infty)$ such that
\begin{equation*}
J\big(\bar{x}(\cdot),\bar{u}(\cdot)\big)
=\inf_{(u(\cdot),x(\cdot))
\in\mathcal{J}_{ad}[0,\infty)}J\big(x(\cdot),
u(\cdot)\big):=V_2(0,x_0).
\end{equation*}

Similarly, we choose $K$ such that \eqref{t32}-\eqref{t34} hold. Let a state feedback control be
$u(\cdot)=Kx(\cdot)+v(\cdot).
$
According to \eqref{t9}, the cost functional \eqref{costf11} can be rewritten as
\begin{equation}\label{costf7}
\begin{aligned}
J_2(v(\cdot))=
\mathbb{E}\int_0^\infty
\begin{pmatrix}
x_1\\
v
\end{pmatrix}^\mathrm{T}
\begin{pmatrix}
\hat{Q}      &S^\mathrm{T}\\
S            &\hat{R}
\end{pmatrix}
\begin{pmatrix}
x_1\\
v
\end{pmatrix} dt,
\end{aligned}
\end{equation}
where
$x_1(\cdot)$ and $\begin{pmatrix}
\hat{Q}      &S^\mathrm{T}\\
S &\hat{R}
\end{pmatrix}$ are given by \eqref{SDE9} and \eqref{t8}, respectively.
Let
\begin{equation}\label{t14}
\begin{pmatrix}
x_1\\
v
\end{pmatrix}=
\begin{pmatrix}
I_h                           &\mathbf{0}\\
-\hat{R}^{-1}S   &I_r
\end{pmatrix}\begin{pmatrix}
x_1\\
v_1
\end{pmatrix}=\begin{pmatrix}
x_1\\
-\hat{R}^{-1}Sx_1+v_1
\end{pmatrix}.
\end{equation}
Then, from \eqref{SDE9}, we have
\begin{subnumcases}{\label{SDE17}}
dx_1(t)=[\tilde{A}_1x_1(t)+B_1v_1(t)]dt
+[\tilde{C}_1x_1(t)
+\tilde{D}_1v_1(t)]dW(t),\label{SDE17a}\\
x_1(0)=x_{10},\label{SDE17b}
\end{subnumcases}
where
$$\tilde{A}_1=A_1-B_1\hat{R}^{-1}
S,\
\tilde{C}_1=C_{11}
-\tilde{D}_1\hat{R}^{-1}S,\
\tilde{D}_1=D_1-C_{12}B_2.
$$
By \eqref{t24} and \eqref{t14}, cost functional \eqref{costf7} becomes
\begin{equation}\label{costf9}
\begin{aligned}
J_3(v_1(\cdot))=
\mathbb{E}\Big[\int_0^\infty
\big(x^\mathrm{T}_1(t)\tilde{Q}x_1(t)
+v_1^\mathrm{T}(t)\hat{R}v_1(t) \big) dt\Big],
\end{aligned}
\end{equation}
where $\tilde{Q}=\hat{Q}-S^\mathrm{T}\hat{R}^{-1}S$.
Define
\begin{equation*}
\begin{aligned}
\mathcal{U}_{loc}[0,\infty)=\bigcap_{T>0}
\mathcal{U}_{ad}[0,T],
\end{aligned}
\end{equation*}
where $\mathcal{U}_{ad}[0,T]$ is defined by \eqref{admiss2}.
Next, we introduce the admissible control set
\begin{equation*}
\begin{aligned}
\mathcal{U}_{ad}[0,\infty)=&\Big\{
u(\cdot)\in\mathcal{U}_{loc}[0,\infty)
\ \Big\vert\ \mathbb{E}\int_0^\infty\vert u(t)\vert ^2dt<\infty,\ \text{the corresponding solution to}\\ &\text{equation \eqref{SDE17} satisfies} \ \mathbb{E}\int_0^\infty\vert x(t)\vert ^2dt<\infty \Big\}.
\end{aligned}
\end{equation*}

\textbf{Problem $\rm (IFP)_K$}. For given $x_{10}\in\mathbb{R}^h$, find optimal control $\bar{v}_1(\cdot)\in \mathcal{U}_{ad}[0,\infty)$ such that
$$J_3(\bar{v}_1(\cdot))=\inf_
{v_1(\cdot)
\in\mathcal{U}_{ad}[0,\infty)}
J_3(v_1(\cdot)):=V_3(0,x_{10}).
$$
Similar to Lemma \ref{lem1}, when $K$ is chosen, we can proof that Problem (IFSP) is equivalent to Problem ${\rm (IFP)_K}$. In the following, we will solve Problem ${\rm (IFP)_K}$ which is a stochastic normal LQ control problem with the order $h< n$.
Let us introduce the following notations
\begin{equation*}
\begin{aligned}
&\mathcal{X}[0,T]=\big\{x(\cdot)\ \vert \ x(\cdot)\ \text{is an} \ \mathcal{F}_t\text{-adapted\ process\ such\ that}\ \mathbb{E}[\int_0^T\vert x(t)\vert^2dt]<\infty \big\},\\
&\mathcal{X}_{loc}[0,\infty)=
\bigcap_{T>0}\mathcal{X}[0,T],\\
&\mathcal{X}[0,\infty)=\big\{x(\cdot)
\in\mathcal{X}_{loc}[0,\infty)\ \vert \ \mathbb{E}[\int_0^\infty\vert x(t)\vert^2dt]<\infty\big\}.
\end{aligned}
\end{equation*}
Note that for any $v_1(\cdot)\in\mathcal{U}_{loc}[0,\infty)$ satisfying $\mathbb{E}[\int_0^\infty\vert v_1(t)\vert ^2dt]<\infty$, the solution $x_1(\cdot)$ to equation \eqref{SDE17} might be just in $\mathcal{X}_{loc}[0,\infty)$. In other words, the corresponding trajectory $x_1(\cdot)$ may not satisfy $\mathbb{E}[\int_0^\infty\vert x_1(t)\vert ^2dt]<\infty$. Thus, the above infinite-horizon cost functional $J_3(v_1(\cdot))$ might not be well-defined.
To guarantee the well-posedness of Problem ${\rm (IFP)_K}$,
we give some notions of controllability introduced by \cite{peng1994backward}. These notions measure the ability of the stochastic normal system \eqref{SDE17} that steers a trajectory from a given initial point $x_{10}$ to a given terminal point
\begin{equation}\label{t15}
x_1(T)=\xi\in L^2_{\mathcal{F}_T}(\Omega;\mathbb{R}^h).
\end{equation}
\begin{definition}{(Definition 3.1, \cite{peng1994backward})}
A stochastic control system \eqref{SDE17a} is called exactly terminal-controllable on $[0,T]$ if for any $\xi\in  L^2_{\mathcal{F}_T}(\Omega;\mathbb{R}^h)$, there exists at least one admissible control $v_1(\cdot)\in \mathcal{U}_{ad}[0,T]$ such that the corresponding trajectory $x_1(\cdot)$ satisfies the terminal condition \eqref{t15}.
\end{definition}
\begin{definition}{(Definition 3.2, \cite{peng1994backward})}
A stochastic control system \eqref{SDE17a} is called exactly controllable on $[0,T]$ if for any $x_{10}\in\mathbb{R}^h$ and $\xi\in  L^2_{\mathcal{F}_T}(\Omega;\mathbb{R}^h)$, there exists at least one admissible control $v_1(\cdot)\in \mathcal{U}_{ad}[0,T]$ such that the corresponding trajectory $x_1(\cdot)$ satisfies the initial condition \eqref{SDE17b} as well as the terminal condition \eqref{t15}.
\end{definition}
\begin{lemma}{(Theorem 3.1,
\cite{peng1994backward})}\label{lem2}
For any $T$, the system \eqref{SDE17a} is exactly terminal-controllable on $[0,T]$ if and only if $rank\ ( \tilde{D}_1)=h.
$
In this case, we can use a simple linear transformation
$v_1=M_2\begin{pmatrix}
z\\
w
\end{pmatrix}+K_1x_1
$ to convert \eqref{SDE17a} into an equivalent form:
\begin{equation}\label{BSDE}
-dx_1(t)=[Px_1(t)+P_1z(t)+\tilde{B}_1w(t)]
dt-z(t)dW(t).
\end{equation}
\end{lemma}
Here, $M_2$ is an invertible $r\times r$-matrix such that $\tilde{D}_1M_2=\begin{pmatrix}
I_h &\mathbf{0}\end{pmatrix}$, $K_1$ is an $r\times h$-matrix such that $\tilde{C_1}+\tilde{D_1}K_1=\mathbf{0}_h$, $P=-(\tilde{A}_1+B_1K_1)$, $P_1=-B_1M_2\begin{pmatrix}
I_h\\
\mathbf{0}
\end{pmatrix}$,
$\tilde{B}_1=-B_1M_2\begin{pmatrix}
\mathbf{0}\\
I_{r-h}
\end{pmatrix}$. $z(\cdot)$ and $w(\cdot)$ are respectively $\mathbb{R}^h$-valued and $\mathbb{R}^{r-h}$-valued processes.
\begin{remark}
For any $\xi \in L^2_{\mathcal{F}_T}(\Omega;\mathbb{R}^h)$, we consider the following BSDE:
\begin{equation*}
\left\{\begin{aligned}
-dx_1(t)&=[Px_1(t)+P_1z(t)+\tilde{B}_1w(t)]
dt-z(t)dW(t),\\
x_1(T)&=\xi.
\end{aligned}
\right.
\end{equation*}
According to Theorem 3.1 in \cite{pardoux1990adapted}, for each $w(\cdot)\in L^2_{\mathcal{F}_t}(0,T;\mathbb{R}^{r-h})$,
there exist a unique pair $(x_1(\cdot),z(\cdot))\in L^2_{\mathcal{F}_t}(0,T;\mathbb{R}^{2h})$ that solves the above BSDE. We denote this pair by $(x_1^w(\cdot),z^w(\cdot))$.
Therefore, the system \eqref{SDE17a} is exactly controllable if and only if for any $x_{10}\in\mathbb{R}^h$ and $\xi \in L^2_{\mathcal{F}_T}(\Omega;\mathbb{R}^h)$, there exists at least one $w(\cdot)\in L^2_{\mathcal{F}_t}(0,T;\mathbb{R}^{r-h})$ such that the solution to the above BSDE satisfies $x_1^w(0)=x_{10}.$
\end{remark}

Once \eqref{SDE17a} is written in the form of \eqref{BSDE}, we have the following algebraic criterion.
\begin{proposition}
{(Theorem 3.2, \cite{peng1994backward})}\label{prop1} For any $T$, the system \eqref{BSDE} is exactly controllable on $[0,T]$ if and only if
$rank \begin{pmatrix}
\tilde{B}_1 &P\tilde{B}_1 &P_1\tilde{B}_1 &PP_1\tilde{B}_1 &P_1P\tilde{B}_1 &\cdots
\end{pmatrix}=h.
$
\end{proposition}

Similar to the deterministic case, we give the following PBH rank criterion that can be used to judge the controllability of the system directly from coefficients.
\begin{proposition}\label{prop2}
For any $T$, the system \eqref{BSDE} is exactly controllable on $[0,T]$ if and only if
\begin{equation}\label{t16}
rank \begin{pmatrix}
sI_h-P &s_1I_h-P_1 &\tilde{B}_1
\end{pmatrix}=h,\ \forall \ s,s_1\in\mathbb{C}.
\end{equation}
\end{proposition}
\noindent{\bf Proof}\quad
Firstly, we prove the necessity part of the statement by contradiction.
We assume that there exist $s$ and $s_1\in \mathbb{C}$ such that rank
$\begin{pmatrix}
sI_h-P & s_1I_h-P_1 &\tilde{B}_1
\end{pmatrix}<h$. Then, there will be a non-zero vector $\beta\in \mathbb{R}^h$ such that
$
\beta^\mathrm{T}\begin{pmatrix}
sI_h-P &s_1I_h-P_1 &\tilde{B}_1
\end{pmatrix}=0,
$
which yields that \begin{equation}\label{t17}
s\beta^\mathrm{T}=\beta^\mathrm{T}P,\ s_1\beta^\mathrm{T}=\beta^\mathrm{T}P_1,\
\beta^\mathrm{T}\tilde{B}_1=0.
\end{equation}
Then we have
\begin{equation}\label{t18}
\begin{aligned}
&\beta^\mathrm{T}P\tilde{B}_1=s\beta^\mathrm{T}
\tilde{B}_1=0,\
\beta^\mathrm{T}P_1\tilde{B}_1=
s_1\beta^\mathrm{T}\tilde{B}_1=0,\\
&\beta^\mathrm{T}PP_1\tilde{B}_1=s
\beta^\mathrm{T}P_1\tilde{B}_1=0,\
\beta^\mathrm{T}P_1P\tilde{B}_1=
s_1\beta^\mathrm{T}P\tilde{B}_1=0,\ \cdots.
\end{aligned}
\end{equation}
From \eqref{t17} and \eqref{t18}, we obtain
$\beta^\mathrm{T}\begin{pmatrix}
\tilde{B}_1 &P\tilde{B}_1 &P_1\tilde{B}_1 &PP_1\tilde{B}_1 &P_1P\tilde{B}_1 &\cdots
\end{pmatrix}=0,
$
which implies that\\
rank $\begin{pmatrix}
\tilde{B}_1 &P\tilde{B}_1 &P_1\tilde{B}_1 &PP_1\tilde{B}_1 &P_1P\tilde{B}_1 &\cdots
\end{pmatrix}<h$. According to Proposition \ref{prop1}, it is contradict with the exact controllability of the system \eqref{BSDE}.
Thus, the necessity of \eqref{t16} holds.

In the next step, we show that the system \eqref{BSDE} is exactly controllable if \eqref{t16} holds. Now, we suppose that the system is not exactly controllable. Consequently,
there exists a $\tilde{x}_{10}\in \mathbb{R}^h$ such that $x^w_1(0)\neq \tilde{x}_{10}$ for any $w(\cdot)\in L^2_{\mathcal{F}_t}(0,T;\mathbb{R}^{r-h})$. Denote the space of all controllable states and uncontrollable states by $X_C$ and $X_{NC}$, respectively.
Making a linear transformation
$x_1=T\tilde{x}_1,\ T=\begin{pmatrix}
T_1 &T_2
\end{pmatrix}.$
Here, each column of matrix $T_1$ belongs to $X_C$ and constitutes the basis of $X_C$, and the matrix $T_2$ is any non-zero matrix that ensures $T$ to be an invertible matrix.
Let $$T^{-1}=\begin{pmatrix}
F_1^\mathrm{T}\\
F_2^\mathrm{T}
\end{pmatrix}.$$
By $T^{-1}T=I_h$, we get $F^\mathrm{T}_2T_1=\mathbf{0},$ which indicates that each column of the matrix $F_2$ belongs to $X_{NC}$.
By Lemma 3.1 in \cite{peng1994backward}, we obtain
$$\Big\{x^w_1(0):w(\cdot)\in L^2_{\mathcal{F}_t}(0,T;
\mathbb{R}^{r-h})\Big\}= span \begin{pmatrix}
\tilde{B}_1 &P\tilde{B}_1 &P_1\tilde{B}_1
&PP_1\tilde{B}_1 &P_1P\tilde{B}_1 &\cdots
\end{pmatrix}.
$$
And then, we derive that each column of $PT_1$ and $P_1T_1$ is in $X_C$, respectively, which yields that $F_2^\mathrm{T}PT_1=0,\ F_2^\mathrm{T}P_1T_1=0.$
It is clear that
\begin{equation*}
\begin{aligned}
&\hat{P}:=T^{-1}PT=\begin{pmatrix}
F_1^\mathrm{T}PT_1  &F_1^\mathrm{T}PT_2\\
\mathbf{0}          &F_2^\mathrm{T}PT_2
\end{pmatrix},\
\hat{P}_1:=T^{-1}P_1T=\begin{pmatrix}
F_1^\mathrm{T}P_1T_1  &F_1^\mathrm{T}P_1T_2\\
\mathbf{0}            &F_2^\mathrm{T}P_1T_2
\end{pmatrix},\\
&\hat{B}:=T^{-1}\tilde{B}_1=\begin{pmatrix}
F_1^\mathrm{T}\tilde{B}_1\\
\mathbf{0}
\end{pmatrix},
\end{aligned}
\end{equation*}
then, we have
\begin{equation*}
\begin{aligned}
&\text{rank} \begin{pmatrix}
sI_h-\hat{P} &s_1I_h-\hat{P}_1 &\hat{B}
\end{pmatrix}\\
=& \text{rank}
\begin{pmatrix}
sI-F_1^\mathrm{T}PT_1 &-F_1^\mathrm{T}PT_2 &s_1I-F_1^\mathrm{T}P_1T_1 &-F_1^\mathrm{T}P_1T_2 &F_1^\mathrm{T}\tilde{B}_1 \\
\mathbf{0} &sI-F_2^\mathrm{T}PT_2
&\mathbf{0} &s_1I-F_2^\mathrm{T}P_1T_2
&\mathbf{0}
\end{pmatrix}.
\end{aligned}
\end{equation*}
It is noted that when $s$ takes the eigenvalue of $F_2^\mathrm{T}PT_2$ and $s_1$ takes the eigenvalue of $F_2^\mathrm{T}P_1T_2$,\\
rank $\begin{pmatrix}
sI_h-\hat{P} &s_1I_h-\hat{P}_1 &\hat{B}
\end{pmatrix}<h$, which contradicts to the following fact:
\begin{equation*}
\begin{aligned}
\text{rank} \begin{pmatrix}
sI_h-\hat{P} &s_1I_h-\hat{P}_1 &\hat{B}
\end{pmatrix}=
&\text{rank}\ T^{-1}
\begin{pmatrix}
sI_h-P &s_1I_h-P_1 &\tilde{B}_1
\end{pmatrix}
\begin{pmatrix}
T &0 &0\\
0 &T &0\\
0 &0 &I_r
\end{pmatrix}\\
=&\text{rank} \ \begin{pmatrix}
sI_h-P &s_1I_h-P_1 &\tilde{B}_1
\end{pmatrix}=h, \ \forall \ s,s_1 \in \mathbb{C}.
\end{aligned}
\end{equation*}
Therefore, we conclude the conclusion.
\hfill$\Box$

In the next step, we present the optimal feedback control of Problem (IFSP) via the following algebraic Riccati equation (ARE)
\begin{equation}\label{ARE3}
\begin{aligned}
&P_0\tilde{A}_1+\tilde{A}_1^{\mathrm{T}}P_0
+\tilde{C}_{1}^
{\mathrm{T}}P_0\tilde{C}_{1}+\tilde{Q}\\
&-(P_0B_1+{\tilde{C}_{1}}^{\mathrm{T}}P_0\tilde{D}_1)
(\hat{R}+\tilde{D}_1^{\mathrm{T}}P_0\tilde{D}_1)^{-1}
(B_1^{\mathrm{T}}P_0+\tilde{D}_1
^{\mathrm{T}}P_0\tilde{C}_{1})=0.
\end{aligned}
\end{equation}
It is noted $P_0$ is the limit of the Riccati equation $P(\cdot;T)\in C(0,T;\mathbb{S}_+^h)$ (see \eqref{DRE1}) deduced from the normal LQ control problem with finite horizon. Therefore, $P_0 \in C(0,T;\mathbb{S}^h_+)$.
Before going further, let us introduce the following assumptions.
\begin{itemize}
\item (H2) $Q>0, \ R>0.$

\item (H2)$'$ $
Q\geq0, \ \text{rank} \left(I_r+KN_1\begin{pmatrix}
\mathbf{0}_{h\times r}\\
-B_2
\end{pmatrix}\right)=r,\ R>0.
$
\item (H3) For fixed $T$, the system \eqref{SDE17} is exactly controllable on $[0,T]$.
\end{itemize}


Similar to (H1) (or (H1)$'$), under (H2) (or (H2)$'$), we have $\tilde{Q}\geq0,\ \hat{R}>0$. Therefore, analogous to Theorem \ref{theorem2}, we obtain the following theorem.
\begin{theorem}\label{theorem5}
Under {\rm (H2)} (or {\rm (H2)$'$}) and {\rm (H3)}, once $K$ is chosen such that \eqref{t32}-\eqref{t34} hold, Problem (IFSP) has a unique optimal control which has the following feedback form
\begin{equation*}\label{optimalv4}
\bar{u}(t)=\big[K+(\Lambda_0-\hat{R}^{-1}S)
\begin{pmatrix}
I_h     &\mathbf{0}_{h\times (n-h)}
\end{pmatrix}N_1^{-1}\big]\bar{x}(t),
\end{equation*}
where
$\Lambda_0=-(\hat{R}+\tilde{D}_1
^{\mathrm{T}}P_0\tilde{D}_1)^{-1}
(B_1^{\mathrm{T}}P_0+\tilde{D}_1
^{\mathrm{T}}P_0\tilde{C}_{1}).$
Moreover, the minimum value of the cost functional is given by
\begin{equation*}\label{costf6}
V_2(0,x_0)=x_0^\mathrm{T}M_1^\mathrm{T}
\begin{pmatrix}
I_h\\
\mathbf{0}
\end{pmatrix}P_0
\begin{pmatrix}
I_h &\mathbf{0}
\end{pmatrix}M_1x_0.
\end{equation*}
\end{theorem}
\noindent{\bf Proof}\quad
From Theorem \ref{theorem6}, Problem ${\rm (IFP)_K}$ has a unique optimal control $\bar{v}_1(\cdot)$, and ARE \eqref{ARE3} admits a solution $P_0\in C(0,T;\mathbb{S}^h_+)$.
Furthermore, the optimal control of Problem ${\rm (IFP)_K}$ is given by
\begin{equation*}\label{optimalv5}
\bar{v}_1(\cdot)=\Lambda_0\bar{x}_1(\cdot),
\end{equation*}
and
$V_3(0,x_{10})=x_{10}^\mathrm{T}
P_0x_{10}.
$
According to \eqref{con1}, \eqref{t7} and \eqref{t14},
we derive the optimal control of Problem (IFSP) as follows
\begin{equation*}
\begin{aligned}
\bar{u}(\cdot)=&K\bar{x}(\cdot)+\bar{v}(\cdot)
=K\bar{x}(\cdot)+(
\Lambda_0-\hat{R}^{-1}S)\bar{x}_1(\cdot)\\
=&\big[K+(\Lambda_0-\hat{R}^{-1}S)
\begin{pmatrix}
I_h     &\mathbf{0}_{h\times (n-h)}
\end{pmatrix}N_1^{-1}\big]\bar{x}(\cdot),\quad t\geq 0,
\end{aligned}
\end{equation*}
and the minimum value of the cost functional
$$V_2(0,x_{0})=V_3(0,x_{10})=x_0^\mathrm{T}M_1^\mathrm{T}
\begin{pmatrix}
I_h\\
\mathbf{0}
\end{pmatrix}P_0
\begin{pmatrix}
I_h &\mathbf{0}
\end{pmatrix}M_1x_0.$$
\hfill$\Box$
\begin{remark}
In the stochastic normal LQ control problem, we usually introduce a Riccati equation with the order $n$ which equals to the order of the system. However, for the singular control problem, we only have a Riccati equation with the order $h< n$. This will simplify the calculation when $h$ is much smaller than $n$.
\end{remark}
\section{Example}
In this section, we consider two examples related to stochastic singular systems.
\begin{example}\label{example1}
Consider the stochastic singular system
\begin{equation*}
\left\{\begin{aligned}
\begin{pmatrix}
1       &0    &0\\
2       &1    &0\\
0       &1    &0
\end{pmatrix}dx(t)&=
\left[\begin{pmatrix}
1       &-1    &0\\
2       &0     &2\\
0       &2     &2
\end{pmatrix}x(t)+\begin{pmatrix}
3\\
0\\
-1
\end{pmatrix}u(t)\right]dt
+\left[\begin{pmatrix}
1\\
4\\
2
\end{pmatrix}u(t)\right]dW(t),\\
\begin{pmatrix}
1       &0    &0\\
2       &1    &0\\
0       &1    &0
\end{pmatrix}x(0) &=  \begin{pmatrix}
1 &2 &0
\end{pmatrix}^\mathrm{T},
\end{aligned}
\right.
\end{equation*}
with cost functional
\begin{equation*}
\begin{aligned}
J\big(u(\cdot),x(\cdot)\big)
=\frac{1}{2}
\mathbb{E}\big[\int_0^1\big
(x^\mathrm{T}(t)x(t)
+u^2(t)\big) dt\big]+\frac{1}{2}\mathbb{E}\left[
x^\mathrm{T}(1)\begin{pmatrix}
5       &2    &0\\
2       &2    &0\\
0       &0    &0
\end{pmatrix}x(1)\right].
\end{aligned}
\end{equation*}
In this system, $n=3$ and rank $(E)=2$. It can be verified that $\vert A+\lambda E\vert\equiv0$, for any $\lambda\in\mathbb{C}$. Thus, there do not exist $M$ and $N$ such that \eqref{t20} holds. According to Theorem \ref{C0}, the system itself does not have a unique solution. However, we can choose
$$K=(0,0,1),\ M_1=\frac{1}{5}\begin{pmatrix}
-1       &3    &-3\\
-2       &1    &4\\
2        &-1   &1
\end{pmatrix}, \
N_1=I_3,$$
such that the system satisfies \eqref{t32}-\eqref{t34} and (H1).
Here, $G=0,\ A_1=\begin{pmatrix}
1       &-1    \\
0       &2
\end{pmatrix}$, $C_{11}=\mathbf{0}_2$,
$C_{12}=\begin{pmatrix}
1\\
2
\end{pmatrix}$, $C_{22}=0$, $B_1=
\begin{pmatrix}
0\\
-2
\end{pmatrix}$, $B_2=1$, $D_1=\begin{pmatrix}
1\\
2
\end{pmatrix}$, $D_2=0$, $Q=I_3$, $R=1$, $H=I_3$. And then, we obtain $\hat{R}=1$, $\tilde{D}_1=\mathbf{0}_{
2\times1},\
S=\mathbf{0}_{1\times2}$, $\hat{Q}=I_2$, $\hat{H}=\begin{pmatrix}
5  &2\\
2  &2
\end{pmatrix}$. According to Theorem 4.2, we have $\Psi=
(-2P_{21},-2P_{22})$, where $P:=\begin{pmatrix}
P_{11}  &P_{12}\\
P_{21}  &P_{22}
\end{pmatrix}$ satisfies the Riccati equation \eqref{Riccati eq} with $P_{12}=P_{21}$. Therefore, the optimal control is given by
\begin{equation*}
\bar{u}(\cdot)=
\big[K-\Psi(t)\begin{pmatrix}
1  &0 &0\\
0  &1 &0
\end{pmatrix}
\big]\bar{x}(\cdot).
\end{equation*}
In the following Figure \ref{figure1}, we give the numerical solution of $P_{11}$, $P_{12}(=P_{21})$ and $P_{22}$.
\end{example}
\begin{figure}[H]
\centering
\includegraphics[width=0.4\textwidth]{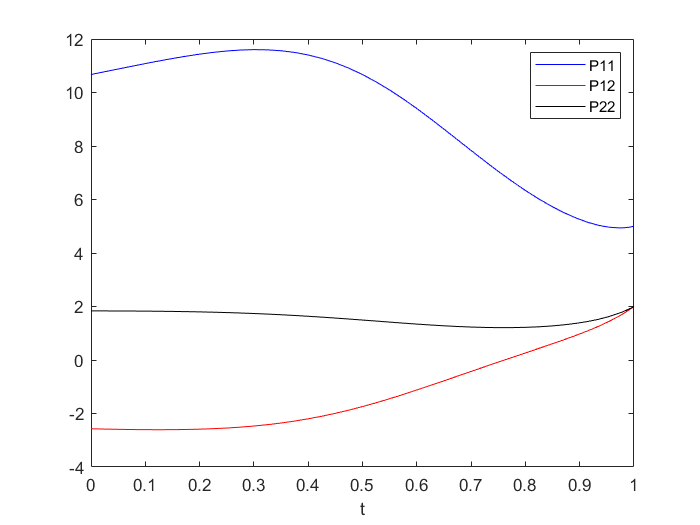}
\caption{The numerical solution of $P$}
\label{figure1}
\end{figure}
\begin{example}\label{example2}
Consider the stochastic singular system
\begin{equation*}
\left\{\begin{aligned}
Edx(t)&=[Ax(t)+Bu(t)]dt
+[Cx(t)+Du(t)]dW(t),\\
Ex(0) &= x_0,\ \begin{pmatrix}
\mathbf{0}_{2\times 1} &I_2
\end{pmatrix}x_0=\mathbf{0}_{2\times1},
\end{aligned}
\right.
\end{equation*}
with the cost functional given by
\begin{equation*}
J\big(x(\cdot),u(\cdot)\big)=
\mathbb{E}\Big[\int_0^\infty
\big(x(t)^\mathrm{T}Qx(t)+ u(t)^\mathrm{T}Ru(t)
\big)dt\Big].
\end{equation*}
where \begin{equation*}
\begin{aligned}
&E=\begin{pmatrix}
1       &0    &0\\
0       &0    &0\\
0       &0    &0
\end{pmatrix}, \
A=\begin{pmatrix}
-1       &0    &0\\
0       &1    &0\\
0       &0    &1
\end{pmatrix},\
B=\begin{pmatrix}
0     &1     &0\\
0     &0     &0\\
0     &0     &1
\end{pmatrix},\\
&C=
D=\begin{pmatrix}
1    &0    &0\\
0    &0    &0\\
0    &0    &0
\end{pmatrix},\
Q=R=I_3.
\end{aligned}
\end{equation*}
On the one hand, let $M_1=N_1=I_3,\ K=\mathbf{0}_{3}$, we can check that \eqref{t32}-\eqref{t34} hold
with
$A_1=-1,\ B_1=\begin{pmatrix}
0 &1 &0
\end{pmatrix},\
B_2=\begin{pmatrix}
0   &0   &0\\
0   &0   &1
\end{pmatrix},\ C_{11}=1,\ C_{12}=C_{21}^\mathrm{T}=
\mathbf{0}_{1\times2},\ C_{22}=\mathbf{0}_{2},\
D_1=\begin{pmatrix}
1   &0   &0
\end{pmatrix},\ D_2=\mathbf{0}_{2\times3}$. From \eqref{t8}, we have:
$$\hat{Q}=1, S=\mathbf{0}_{1\times 3},\ \hat{R}=\begin{pmatrix}
1  &0   &0\\
0  &1   &0\\
0  &0   &2
\end{pmatrix}.$$
By \eqref{SDE17}, we get $\tilde{A}_1=-1,\ \tilde{C_1}=1,\ \tilde{D}_1=D_1,\ \tilde{Q}=1$. Now we apply Lemma \ref{lem2} by choosing $M_2=I_3,\ K_1=\begin{pmatrix}
-1 &0 &0\end{pmatrix}$ to have $P=1,\ P_1=0,\ \tilde{B}_1=\begin{pmatrix}
-1 &0\end{pmatrix}$. Then, rank $\begin{pmatrix}
s-1 &s_1 &-1 &0
\end{pmatrix}=1,\ \forall\ s,\ s_1\in\mathbb{C}$, which yields that (H3) holds by Proposition \ref{prop2}. Thus, we can use Theorem \ref{theorem5} to obtain the optimal feedback control. To do this, by \eqref{ARE3}, we introduce the following ARE:
\begin{equation}\label{ARE5} -P_0-P_0+P_0+1-(P_0B_1+P_0D_1)
(\hat{R}+D_1^\mathrm{T}P_0D_1)^{-1}
(B_1^\mathrm{T}P_0+D_1^\mathrm{T}P_0)
=0,
\end{equation}
which is equivalent to
\begin{equation*}
P_0^3+3P_0^2-1=0.
\end{equation*}
By Cardano formula, it follows that
$
P_{01}=2\cos\frac{2}{9}\pi-1,\
P_{02}=2\cos\frac{4}{9}\pi-1,\
P_{03}=2\cos\frac{8}{9}\pi-1.
$
Noting that we need $P_0\geq0$, it means that the unique solution of \eqref{ARE5} is $P_0=2\cos\frac{2}{9}\pi-1$. Therefore, according to Theorem \ref{theorem5}, the optimal control is given by $\bar{u}=\Lambda_0\begin{pmatrix}
1 &\mathbf{0}_{1\times 2}
\end{pmatrix}\bar{x}$,
where
$
\Lambda_0=-(\hat{R}+D_1^\mathrm{T}
P_0D_1)^{-1}(B_1^\mathrm{T}
P_0+D_1^\mathrm{T}
P_0)
=-\begin{pmatrix}
1-\frac{1}{2\cos\frac{2}{9}\pi}&
2\cos\frac{2}{9}\pi-1&
0
\end{pmatrix}^\mathrm{T}
$. In addition, from \eqref{SDE8}, \eqref{t14} and Theorem \ref{theorem6}, we have:
$$\bar{x}_2=-B_2\bar{v}=-B_2(
-\hat{R}^{-1}
S^\mathrm{T}\bar{x}_1+\bar{v}_1)
=-B_2\Lambda_0 \bar{x}_1.$$ And then, we obtain:
\begin{equation*}
\bar{x}=\begin{pmatrix}
1\\
-B_2\Lambda_0
\end{pmatrix}\bar{x}_1,
\end{equation*}
where $\bar{x}_1$ solves the following SDE from Theorem \ref{theorem6}
\begin{equation*}
\left\{\begin{aligned}
d\bar{x}_1(t)=&(-1
+B_1\Lambda_0)\bar{x}_1(t)dt +(1+D_1\Lambda_0)\bar{x}_1(t)dW(t),\\
\bar{x}_1(0)=&
\begin{pmatrix}
1  &0   &0
\end{pmatrix}x_0:=x_{10}.
\end{aligned}
\right.
\end{equation*}
Thus, we derive $\bar{x}_1(t)=x_{10}\exp^{-1.7450t
+0.6527W(t)}$
and
$\min_{u}J=x_{10}^2P_0=(2\cos
\frac{2}{9}\pi-1)x_{10}^2.
$

On the other hand, the system is equivalent to the following equation:
\begin{equation*}
\left\{\begin{aligned}
dx_1(t)=&(-x_1(t)
+B_1u(t))dt +(x_1(t)+D_1u(t))dW(t),\\
\mathbf{0}_{2\times1}=&x_2+B_2u(t),\\
x_1(0)=& x_{10}.
\end{aligned}
\right.
\end{equation*}
Substituting $x_2=-B_2u$ into the cost functional, we have
\begin{equation*}
J\big(x(\cdot),u(\cdot)\big)
=\mathbb{E}\Big[\int_0^\infty
\big(x_1^2(t)+ u(t)^\mathrm{T}\hat{R}u(t)
\big)dt\Big].
\end{equation*}
By solving the above normal LQ control problem (see Theorem A.1), we can address the original problem.
\end{example}
\section{Conclusions}
This paper works out the necessary and sufficient conditions for ensuring the well-posedness of the stochastic singular system, where both the drift and diffusion terms contain control variables, which is rare in existing literature. Furthermore, we investigate the stochastic singular LQ control problem for both finite and infinite horizons. We transform stochastic singular LQ control problems into stochastic normal LQ control problems using the Kronecker canonical form. To guarantee the finiteness of the infinite-horizon stochastic LQ control problem, we give the PBH rank criterion to verify the exact controllability of the system. Additionally, we derive the feedback optimal control via introducing Riccati equations. Finally, we present solutions for examples illustrating singular LQ control problem.
\section*{Acknowledgments}
This work was supported by the National Key Research and Development Program of China under Grant 2023YFA1009200, the National Natural Science Foundation of China (No. 11831010, 61961160732), the Natural Science Foundation of Shandong Province (No. ZR2019ZD42), the Taishan
Scholars Climbing Program of Shandong (No. TSPD20210302). The work of T. Nie was supported by the National Key R\&D Program of China (No. 2022YFA1006104), the National Natural Science Foundation of China (No. 12022108,11971267,12371447), Natural Science Foundation of Shandong Province  (No. ZR2022JQ01), the Distinguished Young Scholars Program of Shandong University.


\section*{Appendix A}
\setcounter{equation}{0}
\setcounter{subsection}{0}
\renewcommand{\theequation}{A.\arabic{equation}}
\renewcommand{\thesubsection}{A.\arabic{subsection}}
\renewcommand\thetheorem{A.\arabic{theorem}}
\renewcommand\thelemma{A.\arabic{lemma}}
%
%
We consider the following controlled SDE
\begin{equation}\label{SDE28}
\left\{\begin{aligned}
dx_1(t)=&[\tilde{A}_1x_1(t)+B_1v_1(t)]dt
+[\tilde{C}_1x_1(t)
+\tilde{D}_1v_1(t)]dW(t),\ t>0,\\
x_1(0)=&x_{10},
\end{aligned}
\right.
\end{equation}
and the cost functional
\begin{equation}\label{costf10}
\begin{aligned}
J_3(v_1(\cdot))=
\mathbb{E}\Big[\int_0^\infty
\big(x^\mathrm{T}_1(t)\tilde{Q}x_1(t)
+v_1^\mathrm{T}(t)\hat{R}v_1(t) \big) dt\Big].
\end{aligned}
\end{equation}

\textbf{Problem (IFP)}. Find the optimal control $\bar{v}_1(\cdot)\in \mathcal{U}_{ad}[0,\infty)$ such that
$$J_3(\bar{v}_1(\cdot))=\inf_
{v_1(\cdot)
\in\mathcal{U}_{ad}[0,\infty)}
J_3(v_1(\cdot)):=V_3(0,x_{10}).
$$
We introduce the following ARE
\begin{equation}\label{ARE4}
\begin{aligned}
&P_0\tilde{A}_1+\tilde{A}_1^{\mathrm{T}}P_0
+\tilde{C}_{1}^
{\mathrm{T}}P_0\tilde{C}_{1}+\tilde{Q}\\
&-(P_0B_1+{\tilde{C}_{1}}^{\mathrm{T}}P_0\tilde{D}_1)
(\hat{R}+\tilde{D}_1^{\mathrm{T}}P_0\tilde{D}_1)^{-1}
(B_1^{\mathrm{T}}P_0+\tilde{D}_1
^{\mathrm{T}}P_0\tilde{C}_{1})=0.
\end{aligned}
\end{equation}
Define
$\Lambda_0=-(\hat{R}+\tilde{D}_1
^{\mathrm{T}}P_0\tilde{D}_1)^{-1}
(B_1^{\mathrm{T}}P_0+\tilde{D}_1
^{\mathrm{T}}P_0\tilde{C}_{1}).$
Then, we have the following theorem.
\begin{theorem}\label{theorem6}
Let {\rm (H2)} (or {\rm (H2)$'$}) and {\rm (H3)} hold. ARE \eqref{ARE4} admits a solution $P_0 \in C(0,T;\mathbb{S}^h_+)$,
and Problem (IFP) has a unique optimal control $\bar{v}_1(\cdot)\in \mathcal{U}_{ad}[0,\infty)$ which satisfies the following feedback form:
\begin{equation*}
\bar{v}_1(t)=\Lambda_0\bar{x}_1(t),\ t\geq 0,
\end{equation*}
where $\Lambda_0=-(\hat{R}+\tilde{D}_1
^{\mathrm{T}}P_0\tilde{D}_1)^{-1}
(B_1^{\mathrm{T}}P_0+\tilde{D}_1
^{\mathrm{T}}P_0C_{11})$.
Moreover, the optimal state $\bar{x}_1(\cdot)$ solves
\begin{equation*}
\left\{\begin{aligned}
&d\bar{x}_1(t) = (\tilde{A}_1
+B_1\Lambda_0)\bar{x}_1(t)dt
+(\tilde{C}_1
+\tilde{D}_1\Lambda_0)\bar{x}_1(t)dW(t),\\
&\bar{x}_1(0)=x_{10},
\end{aligned}
\right.
\end{equation*}
and the minimal cost functional is given by
\begin{equation*}
V_3(0,x_{10})=
x_{10}^\mathrm{T}P_0x_{10}.
\end{equation*}
\end{theorem}
\noindent{\bf Proof}\quad
The proof is similar to \cite{huang2015linear}. From (H3), we know that for any $x_{10}$, there exists at least one control $\hat{v}_1(\cdot)$ such that the corresponding trajectory $\hat{x}_1(\cdot)$ of SDE \eqref{SDE28} satisfies $\hat{x}_1(T)=0$.
Then we define an admissible control as follows
\begin{equation*}
\tilde{v}_1(t)=\left\{
\begin{aligned}
&\hat{v}_1(t),\ t\in [0,T],\\
&0, \ \quad \ \ t\in (T,\infty),
\end{aligned}
\right.
\end{equation*}
and the corresponding trajectory is given by
\begin{equation*}
\tilde{x}_1(t)=\left\{
\begin{aligned}
&\hat{x}_1(t),\ t\in [0,T],\\
&0, \ \quad \ \ t\in (T,\infty).
\end{aligned}
\right.
\end{equation*}
Then, we obtain
$$J_3(\tilde{v}_1(\cdot))=\mathbb{E}\int_0
^T\big[\hat{x}_1(t)^\mathrm{T}\tilde{Q}\hat{x}_1(t)+
\hat{v}_1(t)^\mathrm{T}\hat{R}\hat{v}_1(t)\big]dt <\infty.
$$
Therefore, the admissible control set $\mathcal{U}_{ad}[0,\infty)$ is nonempty and the value function $V_3(0,x)$ is bounded by $J_3(\tilde{v}_1(\cdot))$.
Meanwhile, $J_3(v_1(\cdot))$ is a strictly convex and coercive functional under (H2) (or (H2)$'$), which yields that Problem {\rm (IFP)} has a unique optimal control. And then, we can verify that the value function $V_3(0,x)$ for a normal LQ control problem is a quadratic function with respect to $x$. As a consequence, there exists a unique bounded $\bar{P}$ such that $V_3(0,x_{10})=x_{10}^\mathrm{T}\bar{P}x_{10}$. We now show that $\bar{P}$ satisfies equation \eqref{ARE4}. Firstly, we consider the finite-horizon LQ control problem, which is to search
an optimal control $\bar{v}(\cdot)\in \mathcal{U}_{ad}[0,T]$ such that
$$
J(\bar{v}(\cdot);T)
=\inf_{v(\cdot)\in \mathcal{U}_{ad}[0,T]} J(v(\cdot);T):=V(0,x_{10};T),$$
where
\begin{equation}
J(v(\cdot);T)=\mathbb{E}\big[\int_0^T
\big(x^\mathrm{T}(t;T)\tilde{Q}x(t;T)
+v^\mathrm{T}(t)\hat{R}v(t)
\big)dt\big],
\end{equation}
and $x(\cdot;T)$ satisfies equation \eqref{SDE28} driven by $v(\cdot)$.
We denote the above control problem by Problem ${\rm (FP)_T}$. From Theorem 6.1 in \cite{yong1999stochastic}, under (H2) (or (H2)$'$), there is a unique optimal control $\bar{v}(\cdot)\in \mathcal{U}_{ad}[0,T]$ such that
\begin{equation}
\inf_{v(\cdot)\in\mathcal{U}_{ad}[0,T]}
J(v(\cdot);T)=
J(\bar{v}(\cdot);T)=x_{10}^\mathrm{T}
P(0;T)x_{10}, \ \forall \ x_{10}\in\mathbb{R}^h,
\end{equation}
where $P(\cdot;T)\in C(0,T;\mathbb{S}^h_+)$ satisfies the following Riccati equation
\begin{equation}\label{DRE1}
\left\{\begin{aligned}
&\dot{P}(s;T)+P(s;T)\tilde{A}_1+
\tilde{A}_1^\mathrm{T}P(s;T)+
\tilde{C}_1^\mathrm{T}P(s;T)\tilde{C}_1
+\tilde{Q}-[P(s;T)B_1+\tilde{C}_1^\mathrm{T} P(s;T)\tilde{D}_1]\\
&\cdot[\hat{R}+\tilde{D}_1^\mathrm{T}P(s;T)\tilde{D}_1]^{-1}
[B_1^\mathrm{T}P(s;T)+\tilde{D}_1^\mathrm{T}P(s;T)\tilde{C}_1]=0,
\quad s\in [0,T],\\
&P(T;T)=0.
\end{aligned}
\right.
\end{equation}
Moreover, the optimal control $\bar{v}(\cdot)$ is given by
\begin{equation*}
\bar{v}(\cdot)=\Lambda(\cdot;T)
\bar{x}(\cdot;T),
\end{equation*}
where $\Lambda(\cdot;T)=-[\hat{R}
+\tilde{D}_1^\mathrm{T}
P(\cdot;T)\tilde{D}_1]^{-1}
[B_1^\mathrm{T}P(\cdot;T)
+\tilde{D}_1^\mathrm{T}
P(\cdot;T)\tilde{C}_1]$.
Due to $\tilde{Q}\geq 0, \hat{R}>0$, for $0\leq T\leq \bar{T}<\infty$, we have:
\begin{equation*}
\begin{aligned}
x_{10}^\mathrm{T}P(0;T)x_{10}
=&V(0,x_{10};T)=
\inf_{u(\cdot)\in \mathcal{U}_{ad}[0,T]} J(u(\cdot);T)
\leq\inf_{u(\cdot)\in \mathcal{U}_{ad}[0,T]} J(u(\cdot);\bar{T})\\
\leq&\inf_{u(\cdot)\in \mathcal{U}_{ad}[0,\bar{T}]} J(u(\cdot);\bar{T})
=V(0,x_{10};\bar{T})
=x_{10}^\mathrm{T}P(0;\bar{T})x_{10}\\
\leq&\inf_{u(\cdot)\in \mathcal{U}_{ad}[0,\infty)} J(u(\cdot);\bar{T})
\leq\inf_{u(\cdot)\in \mathcal{U}_{ad}[0,\infty)} J(u(\cdot))\\
=&V_3(0,x_{10})
=x_{10}^\mathrm{T}\bar{P}x_{10},\ \text{for any} \ x_{10}\in\mathbb{R}^h.
\end{aligned}
\end{equation*}
Noting that $x_{10}$ can be chosen arbitrarily, we have
$0\leq P(0;T)\leq P(0;\bar{T})\leq \bar{P}$, for $0\leq T\leq \bar{T}<\infty$,
which implies that $\big\{P(0;T)\big\}_{T\geq0}$ is an increasing sequence of $T$ and is bounded by $\bar{P}$. Consequently,  $\lim_{T\rightarrow\infty}P(0;T)$ exists, and we denote
\begin{equation}\label{limT}
P_\infty
:=\lim_{T\rightarrow\infty}P(0;T).
\end{equation}
Define
$P(\cdot,\infty):=
\mathop{\lim}\limits_{T\rightarrow\infty}P(\cdot;T),\
\Lambda(\cdot,\infty):=
\mathop{\lim}\limits_{T\rightarrow\infty}
\Lambda(\cdot;T).$
Meanwhile, we introduce the following Riccati equation $\tilde{P}(\cdot)$ on $[0,\infty)$
\begin{equation}\label{ARE2}
\left\{\begin{aligned}
&-\dot{\tilde{P}}(r)+\tilde{P}(r)\tilde{A}_1
+\tilde{A}_1^\mathrm{T}\tilde{P}(r)+
\tilde{C}_1^\mathrm{T}\tilde{P}(r)\tilde{C}_1
+\tilde{Q}-[\tilde{P}(r)B_1
+\tilde{C}_1^\mathrm{T} \tilde{P}(r)\tilde{D}_1]\\
&\cdot[\hat{R}+\tilde{D}_1^\mathrm{T}
\tilde{P}(r)\tilde{D}_1]^{-1}
[B_1^\mathrm{T}\tilde{P}(r)
+\tilde{D}_1^\mathrm{T}\tilde{P}(r)\tilde{C}_1]
=0,
\quad r\in [0,T],\\
&\tilde{P}(0)=0.
\end{aligned}
\right.
\end{equation}
For any $T>0$, let
$
\hat{P}(s;T)=\tilde{P}(T-s)$. One can check that $\hat{P}(\cdot;T)$ solves \eqref{DRE1}. Then from the uniqueness of \eqref{DRE1} (see Theorem 7.2 in \cite{yong1999stochastic}),
we have
$
P(s;T)=\hat{P}(s;T)=\tilde{P}(T-s),\ s\in [0,T].$
Hence,
$P(0;T)=\tilde{P}(T),\ T\geq 0.
$
Recalling \eqref{limT}, we obtain
$
\lim_{T\rightarrow\infty}\tilde{P}(T)
=P_\infty.$ Integrating both sides of \eqref{ARE2} from $T$ to $T+1$, we get:
\begin{equation*}
\begin{aligned}
&\tilde{P}(T)-\tilde{P}(T+1)+
\int_T^{T+1}\Big\{\tilde{P}(r)\tilde{A}_1
+\tilde{A}_1^\mathrm{T}\tilde{P}(r)+
\tilde{C}_1^\mathrm{T}\tilde{P}(r)\tilde{C}_1
+\tilde{Q}-[\tilde{P}(r)B_1
+\tilde{C}_1^\mathrm{T} \tilde{P}(r)\tilde{D}_1]\\
&\cdot[\hat{R}+\tilde{D}_1^\mathrm{T}
\tilde{P}(r)\tilde{D}_1]^{-1}
[B_1^\mathrm{T}\tilde{P}(r)
+\tilde{D}_1^\mathrm{T}\tilde{P}(r)\tilde{C}_1]
\Big\}dr=0.
\end{aligned}
\end{equation*}
Let $T\rightarrow\infty$, we obtain that $P_\infty$ satisfies ARE \eqref{ARE4}, which means that \eqref{ARE4} has a solution $P_0=P_\infty$.
In addition, we find that
\begin{equation*}
\begin{aligned}
&\lim_{T\rightarrow\infty}P(s;T)
=\lim_{T\rightarrow\infty}\tilde{P}(T-s)
=P_\infty,\\
&\lim_{T\rightarrow\infty}\Lambda(s;T)=
-(\hat{R}+\tilde{D}_1^\mathrm{T}P_\infty \tilde{D}_1)^{-1}
(B_1^\mathrm{T}P_\infty+\tilde{D}_1
^\mathrm{T}P_\infty \tilde{C}_1):=\Lambda_\infty,
\end{aligned}
\end{equation*}
which yields that $\lim_{T\rightarrow\infty}P(s;T)$ and $\lim_{T\rightarrow\infty}\Lambda(s;T)$ are time-invariant constant matrices, respectively.
Let $T\rightarrow\infty$, and define
$\check{x}(\cdot):=\lim_{T\rightarrow\infty}\bar{x}(\cdot;T), \ \check{v}(\cdot):= \lim_{T\rightarrow\infty}\bar{v}(\cdot)
=\Lambda_\infty \check{x}(\cdot)$, then, $\check{x}(\cdot)$ satisfies the following SDE:
\begin{equation*}
\left\{\begin{aligned}
d\check{x}(t)=&(\tilde{A}_1+B_1\Lambda_\infty)
\check{x}(t) dt
+(\tilde{C}_1+\tilde{D}_1\Lambda_\infty)\check{x}(t) dW(t), \ t\geq 0,\\
\check{x}(0)=& x_{10}.
\end{aligned}
\right.
\end{equation*}
From the dominated convergence theorem, we have:
\begin{equation*}
x_{10}^\mathrm{T}\bar{P}x_{10}\geq x_{10}^\mathrm{T}P_\infty x_{10}
=\mathbb{E}\Big[\int_0^\infty
\big(\check{x}^\mathrm{T}(t)\tilde{Q}\check{x}(t)
+\check{v}^\mathrm{T}(t)\hat{R}
\check{v}(t)\big)dt\Big]
\geq x_{10}^\mathrm{T}\bar{P}x_{10},\
\text{for any}\ x\in\mathbb{R}^h,
\end{equation*}
which implies that
$\bar{P}\equiv P_\infty$ and $\bar{P}$ satisfies \eqref{ARE4}. Moreover, from the definition of $\Lambda_0$ and $\Lambda_\infty$, it is obvious that   $\Lambda_0=\Lambda_\infty.$
Hence,
$\bar{x}_1(\cdot)=\check{x}(\cdot),\ \bar{v}_1(\cdot)=\check{v}(\cdot).$
This completes the proof of Theorem  \ref{theorem6}.
\hfill$\Box$
\end{document}